\documentclass[journal,final]{IEEEtran}
\usepackage{amsmath,amsfonts}
\usepackage{algorithm}
\usepackage{algpseudocode}
\usepackage{array}
\usepackage{textcomp}
\usepackage{stfloats}
\usepackage{url}
\usepackage{verbatim}
\usepackage{graphicx}
\usepackage{cite}
\usepackage{xcolor}
\usepackage{caption}
\usepackage{multicol}

\usepackage{dirtytalk}
\usepackage{csquotes}
\AtBeginEnvironment{quote}{\itshape}

\usepackage{booktabs}

\usepackage{mathtools}
\newcounter{tableeqn}[table]

\newcounter{tablesubeqn}[tableeqn]

\usepackage{subfigure}

\newcolumntype{C}[1]{>{\centering\arraybackslash}p{#1}}

\usepackage{pgfplots}
\usetikzlibrary{matrix, shapes, shadows, arrows}
\usepackage{svg}  
\usetikzlibrary{plotmarks}
\usetikzlibrary{arrows.meta}
\usepgfplotslibrary{patchplots}

\newtheorem{theorem}{Theorem}
\newtheorem{proposition}{Proposition}
\newtheorem{corollary}{Corollary}
\newtheorem{lemma}{Lemma}
\newtheorem{remark}{Remark}

\pgfplotsset{compat=1.3}

\newcommand{\cz}[1]{\mathcal{Z}_{#1}}
\newcommand{\cnz}[1]{\mathcal{N}_{#1}}

\begin{document}

\title{Asymptotics of Distances Between Sample Covariance Matrices}

\author{Roberto~Pereira,~\IEEEmembership{Student~Member,~IEEE,}
        Xavier~Mestre,~\IEEEmembership{Senior~Member,~IEEE,}
        and~David~Gregoratti,~\IEEEmembership{Senior~Member,~IEEE}
        
\thanks{R. Pereira and X. Mestre are with the Centre Tecnol\`ogic de Telecomunicacions de Catalunya (CTTC), 
Barcelona (Spain), 
\{roberto.pereira, xavier.mestre\}@cttc.cat}
\thanks{D. Gregoratti is with Software Radio Systems, Barcelona, Spain, david.gregoratti@srs.io}

\thanks{Part of this manuscript has been presented to IEEE GLOBECOM 2021 and IEEE ICASSP 2022.}

\thanks{This work has been partially supported by the European Commission under the Windmill project (contract 813999) and the Spanish Grant PID2021-128373OB-I00 funded by MCIN/AEI/10.13039/501100011033 and by “ERDF A way of making Europe”.}}

\markboth{IEEE Transactions on Signal Processing}%
{Shell \MakeLowercase{\textit{et al.}}: A Sample Article Using IEEEtran.cls for IEEE Journals}


\maketitle

\begin{abstract}

This work considers the asymptotic behavior of the distance between two sample covariance matrices (SCM). A general result is provided for a class of functionals that can be expressed as sums of traces of functions that are separately applied to each covariance matrix. In particular, this class includes very conventional metrics, such as the Euclidean distance or Jeffrery's divergence, as well as a number of other more sophisticated distances recently derived from Riemannian geometry considerations, such as the log-Euclidean metric. In particular, we analyze the asymptotic behavior of this class of functionals by establishing a central limit theorem that allows us to describe their asymptotic statistical law. In order to account for the fact that the sample sizes of two SCMs are of the same order of magnitude as their observation dimension, results are provided by assuming that these parameters grow to infinity while their quotients converge to fixed quantities. Numerical results illustrate how this type of result can be used in order to predict the performance of these metrics in practical machine learning algorithms, such as clustering of SCMs.

\end{abstract}

\begin{IEEEkeywords}
Statistical Analysis, Covariance Matrix  Distance, Random Matrix Theory, Riemannian Geometry.
\end{IEEEkeywords}

\section{Introduction}

A large number of applications in signal processing and machine learning explicitly rely  on a proper quantification of distances among distinct groups of observed data. One of the most important applications is data clustering, which basically identifies groups of affine observations based on some measure of proximity or similarity. Quite often the data to be analyzed is multivariate, in the sense that each observation consists of readings from multiple entities or sensors. Furthermore,  it is usually the case that clustering needs to be performed according to the pattern of dependence between  these multiple readings, rather than proximity between the different observations as such. In other words, the relevant information is contained in the covariance pattern of multivariate observations, rather than the actual measurements. For example, the covariance of the signal that is received by a number of spatially distributed antennas is directly related to the spatial distribution of the corresponding sources, just in the same way that the covariance of a time series is directly connected to the corresponding spectral density. 

Clustering according to covariance matrices is now common in multivariate and functional data analysis \cite{Ieva16, Lee15}. The main application in this setting consists in grouping segments of data that present the same spatial correlation structure. For example,  one can imagine a set of sensors collecting signals in different states,
where discovering and clustering the different covariance of these states is a previous step in order to understand trends, detect changes/failure and generally interpret the high-dimensional data \cite{Hallac17}. 
Another relevant example can be found in multi-user MIMO communications, where it is usually convenient to group users that approximately share the same spatial characteristics 
\cite{adhikary2013joint_spatial_first,xu2014user,sun2017agglomerative,nam2014_jsdm,Czink12,mao2018rate_splitting}. Here, the main objective is to group a number of channel matrices so that channels that are seen from the same spatial locations belong to the same group. Since spatial distribution is directly related to the inter-antenna signal covariance, one can alternatively cluster  these MIMO channels according to their receive covariance matrices. 

Distances between covariance matrices have their own importance beyond the clustering problem. For example, image set classification is largely based on discriminant analysis  on the intra-set covariance matrices \cite{Wang12}. By identifying each image set with its natural second-order statistic, the classification problem can be formulated as discriminating points in the Riemann manifold of positive semi definite matrices \cite{horev2016geometry}.  A similar approach can be followed in diffusion tensor imaging \cite{dryden09}, where the main descriptor is a vortex-depending  covariance matrix; brain-computer interfacing \cite{barachant2011multiclass}, where  covariance  matrices  are used to represent the spatial
information  embedded  in  EEG  signal; or radar/sonar signal processing \cite{Barbaresco08,wong2017mean_euclidean},
where the spatial covariance matrix is used to capture the spatial characteristics of the clutter.
In all these settings, the conventional Euclidean metric is not appropriate for measuring proximity
between the observed covariance matrices, which belong to the set of positive semidefinite
matrices. 
Hence, a number of studies propose to rely on metrics that consider the topological structure of underlying  manifolds~\cite{dryden09,thanwerdas2023n}, e.g.,  Riemannian based distances~\cite{moakher2006symmetric_kl, wong2017mean_euclidean} instead of the classical Euclidean distance.


Unfortunately, in many real-world applications it may become challenging to obtain perfect knowledge of the second order statistics of the underlying processes. Consequently, the inherent distances must be estimated from the corresponding data. This has been the focus of multiple contributions of geometry-aware  metrics for covariance matrices.  The main difficulty here is the fact that the number of observations from which the  covariance matrices need to be estimated is typically low compared to the observation dimension. This makes it extremely difficult to obtain a reliable estimator of the distances in practical settings. 

One possible way of addressing this estimation is by assuming some  structure on the actual covariance matrices. This was the case for instance in~\cite{elvander2018interpolation}, which considered the estimation of a distance measure between Toeplitz covariance matrices based on an optimal mass transport problem. The distance between matrices is then formulated as a
convex optimization problem. 
A similar idea has also been explored in~\cite{ferrante2011maximum}, which proposes using a maximum entropy approach to enhance  a family of high-resolution spectral estimators.
For the more general case where the covariance matrices do not have a specific structure,
some recent contributions \cite{Tiomoko19, tiomoko2019random,couillet2019random,pmlr-v97-tiomoko19a}  have considered the use of random matrix theory results to obtain estimators of the distance that are consistent even when the observation dimension increases with
the sample size. Unfortunately, these works have mainly focused on a specific type of  distance that can be formulated as linear spectral statistics of the product of one covariance
times the inverse of the other. Furthermore, these estimators have only been formulated
for the oversampled regime, where the number of samples available to estimate each covariance
is larger than the observation dimension.
In practice, one usually needs to deal with the undersampled situation, where the number of samples that are 
available to estimate (at least one of) the two covariance matrices is smaller than the 
corresponding observation dimension. To the best of our knowledge, this regime has only been addressed in \cite{Timoko19} by using a polynomial approximation of the functions involved in the distance definition. 

In order to effectively address the undersampled regime, we consider in this paper  distances between \emph{sample} covariance matrices rather than attempting to estimate the distance between the true ones. Our motivation is twofold: on the one hand, distances between sample covariance matrices (SCMs) are relatively easy to implement; on the other hand, these metrics can also be employed in the undersampled regime, a regime for which no consistent estimators generally exist. Our objective is to provide a complete asymptotic statistical characterization of a general class of distances between SCMs, assuming that both the sample size and the observation dimension grow large at the same rate. We will generally consider both the undersampled and the oversampled regime, in the sense that the SCMs will not be required to be positive definite. The asymptotic results will be used to assess the effectiveness of these  measures in practical problems, providing a useful tool to establish which distance performs statistically better for a given scenario. 
It should be pointed out here that it is not the objective of the paper to come up with new consistent estimators for the corresponding distances between the true covariance matrices, which have an interest of their own and are left for future work (see \cite{pereira_icassp23} for a first step in this direction).

We can summarize the main contributions of the present paper as follows. We study the asymptotic behavior of a general class of functions between two SCMs, which can be particularized to several distances and divergences, in both the oversampled and undersampled regime. We derive a new Central Limit Theorem for this type of functionals and particularize the result to several relevant distances between SCMs. We then illustrate the relevance of the derived results by addressing the performance of these metrics in a practical illustrative clustering application. Part of the results presented here were published without proof in a couple of conference
contributions in \cite{roberto21_globecom} and~\cite{roberto22_icassp}. More specifically, \cite{roberto21_globecom} considered a MIMO wireless channel clustering application and focused on the asymptotic regime of the Projection-Frobenius (PF) squared distance between subspaces spanned by MIMO channel matrices. These results were later refined in  \cite{roberto22_icassp}, which considered the consistent estimation of the asymptotic mean and  variance.  In the current work, we provide a generalization of these contributions by considering a family of functionals that includes the aforementioned distances. Furthermore, both complex- and real-valued observations are considered here, and extended proofs are provided.

\section{Statistical Model of the Observations}
\label{sec:statistical_model_observations}
We consider two sets of multidimensional observations of
dimensionality $M$, which are denoted $\mathbf{y}_{1}(n)\in\mathbb{C}
^{M\times1}$ and $\mathbf{y}_{2}(n)\in\mathbb{C}^{M\times1}$ respectively,
$n\in\mathbb{N}$. We assume that the first sample set contains $N_{1}$
observations, whereas the second one is composed of $N_{2}$ observations. We
will denote by $\mathbf{Y}_{1}$and $\mathbf{Y}_{2}$ two matrices of dimensions
$M\times N_{1}$ and $M\times N_{2}$ respectively, which contain the
observations associated to each of these observation sets as columns, that is
\[
\mathbf{Y}_{j}=\left[
\begin{array}
[c]{ccc}
\mathbf{y}_{j}(1) & \ldots & \mathbf{y}_{j}(N_{j})
\end{array}
\right]
\]
for $j\in\{1,2\}$. When the observations are zero mean, the covariance matrix
of these observations are typically estimated using the SCMs, which are defined as
\[
\mathbf{\hat{R}}_{j}=\frac{1}{N_{j}}\mathbf{Y}_{j}\mathbf{Y}_{j}^{H}.
\]
We are interested in measuring the distance between the two SCMs, namely $\mathbf{\hat{R}}_{1}$, $\mathbf{\hat{R}}_{2}$. In
this paper, we will consider the family of quantities that can be mathematically expressed as
\begin{equation}
\hat{d}_{M}(1,2)=\sum_{l=1}^{L}\frac{1}{M}\mathrm{tr}\left[  f_{1}^{(l)}\left(
\mathbf{\hat{R}}_{1}\right)  f_{2}^{(l)}\left(  \mathbf{\hat{R}}_{2}\right)
\right]  \label{eq:definitionhatd}
\end{equation}
for certain functions $f_{1}^{(l)}, f_{2}^{(l)}: \mathbb{C}^{M\times M}\rightarrow\mathbb{C}^{M\times M}$, $l=1,\ldots,L$. 
Here, $L$ is a fixed integer that denotes the number of terms in the definition of $\hat{d}_M(1,2)$. Typically, these functions are understood to be scalar analytic functions applied to the eigenvalues of the Hermitian matrices $\hat{\mathbf{R}}_j$, $j \in \{1,2\}$. Let us particularize this definition to some meaningful choices for the considered distances.

Obviously, any (squared) distance\footnote{Many of the quantities presented here represent the square-form of a distance. Throughout this paper, we will omit this in the notation.} that can be expressed in the form 
\begin{equation}
\hat{d}_{M}(1,2)=\frac{1}{M}\mathrm{tr}\left[  \left( f_{1}(\mathbf{\hat{R}}_{1})  - f_{2}(  \mathbf{\hat{R}}_{2}) \right)^2 \right] \label{eq:d12simple}
\end{equation}
for two matrix-valued functions $f_{1},f_{2}$ can be seen as a particularization of the general functional in (\ref{eq:definitionhatd}). In particular, if  these two functions are chosen so that $f_{j}(\hat{\mathbf{R}}_j)=\hat{\mathbf{R}}_j$, $j\in\{1,2\}$ we will  recover the conventional Euclidean distance
between SCMs, that is
\[
\hat{d}_M^{E}(1,2) = \frac{1}{M}\mathrm{tr}\left[ \left( \mathbf{\hat{R}}_{1}-\mathbf{\hat{R}}_{2} \right)^2 \right]. 
\]
Likewise, the choice $f_{j}(\hat{\mathbf{R}}_j)=\log (\hat{\mathbf{R}}_j)$ will lead to the log-Euclidean distance between SCMs \cite{huang2015log}, 
 whereas the choice $f_{j}(\hat{\mathbf{R}}_j)=(\hat{\mathbf{R}}_j)^\alpha$ for some $\alpha > 0$ will lead to the power-Euclidean distance in \cite{dryden2010power}.


Similarly, after proper normalization,  the symmetrized version of the Kullback-Leibler divergence between two multivariate Gaussians  (usually referred to as Jeffreys divergence \cite{Ali66}) can be expressed by
\begin{equation}
\hat{d}_M^{KL}(1,2) = \frac{1}{2M}\mathrm{tr}\left[\mathbf{\hat{R}}_{1}\mathbf{\hat{R}}_{2}^{-1} \right] + \frac{1}{2M}\mathrm{tr}\left[\mathbf{\hat{R}}_{2}\mathbf{\hat{R}}_{1}^{-1} \right] - 1 \label{eq:defJefferiesdivergence}
\end{equation}
which also conforms to the general expression in (\ref{eq:definitionhatd}). It is worth mentioning that the KL divergence belongs to a larger family of divergences, namely the ($\alpha, \beta$)-divergence~\cite{zhang2004divergence}, that also follows the definition in~(\ref{eq:definitionhatd}). Nonetheless, due to space limitation, we will focus on its particularization to the  divergence above. In fact, we will consider a generalized version of this distance, also valid in the undersampled regime by replacing the inverse $(\cdot)^{-1}$ by the Moore-Penrose pseudoinverse. 
Another alternative would be to replace the functions $f_j^{l}(z) = 1/z$ in the above definition by 
some other function of the form $f_j^{l}(z) = 1/(z + \xi)$ for some $\xi>0$, leading to ``diagonally loaded" versions of the SCMs (with diagonal loading equal to $\xi$). 

Other functionals defined only for positive SCMs accept the same type of generalization to the undersampled regime. For example, one could extend the log-Euclidean metric to the undersampled case by either considering the logarithm of the positive eigenvalue only, or by introducing the diagonally loaded SCMs in the original definition.
To avoid the introduction of yet another parameter ($\xi$) and in the interest of simplicity of presentation, in this paper we will mainly focus on the first type of generalization.

Interestingly enough, the above formulation also particularizes to proximity measures between the column spaces of two sets of observations $\mathbf{Y}_1,\mathbf{Y}_2$ in the  undersampled regime (so that both $\mathbf{Y}_1,\mathbf{Y}_2$ are tall matrices). In  this case, one can measure the similarity between the two subspaces as sum of the squared sines of the principal angles between these subspaces, which can be seen as a squared distance in the Grassmann manifold \cite{zhang2018grassmannian}, i.e.
\[
\hat{d}_M^{SS} (1,2) = \frac{1}{M}\mathrm{tr}\left[ \left( \mathbf{{P}}_{1}-\mathbf{{P}}_{2} \right)^2 \right]
\]
where $\mathbf{P}_i = \mathbf{Y}_i \left( \mathbf{Y}_i^H \mathbf{Y}_i \right)^{-1} \mathbf{Y}_i^H $ is the projection matrix onto the column space of $\mathbf{Y}_i$, $i=1,2$. This squared distance can also be expressed as in (\ref{eq:d12simple}), where now $f_1(z)$ and $f_2(z)$ are functions that such that $f_j(0)=0$ while $f_j(z)=z$ for $z\in\mathbb{R}^+, z>\epsilon$ with $\epsilon>0$ small enough. 
The study of the above quantity has an interest beyond the framework of this paper and can be used to characterize independence tests based on canonical correlation analysis, which typically use $\mathrm{tr}\left[ \mathbf{{P}}_{1}\mathbf{{P}}_{2} \right]$ as the relevant statistic to determine whether the two sets of
observations are statistically  independent (see \cite{han2018unified, bao2019canonical, yang2021limiting} 
for the problem formulation and the asymptotic characterization when the observations are spatially white; results in this paper extend this characterization to the general spatially colored case). 
Similarly, there also exists a wide range of applications in the field of wireless communication that involve the clustering of wireless channels, which can be based on the aforementioned metric~\cite{roberto21_globecom, rs_2rhs, nam2014_jsdm}.

The objective of this paper is to provide an asymptotic characterization of quantities like the ones above when the dimensions of the matrices $M$ and the corresponding sample sizes $N_{1},N_{2}$ increase to infinity at the same rate, so that their quotient converges to a fixed quantity, namely $M/N_{1}\rightarrow c_{1}$, $M/N_{2}\rightarrow c_{2}$ for some $0<c_{1} ,c_{2}<\infty$. The main advantage of this characterization with respect to the more conventional one (which assumes fixed $M$) is the fact that here all the dimensions ($M,N_{1},N_{2}$) are comparable in magnitude even in the asymptotic regime, which makes the analysis more reliable in order to analyze the behavior of $\hat{d}_{M}$ in the finite sample size regime. 
The results provided in this work do not necessarily consider any specific structure on the covariance matrices and are extendable to any metric that falls within the definition described above. Notice that this is different from what was considered previous works found in the literature~\cite{dryden09, barachant2011multiclass, Barbaresco08}, these often aim at studying or designing the best distance tailored to some specific problem.
In what follows, we will carry out the analysis in two steps. First, in Section \ref{sec:firstorder} we will see that under standard statistical assumptions the above distances asymptotically behave as deterministic quantities, which will be referred to as deterministic equivalents.  Then, in Section~\ref{sec:secondorder} we will prove that these distances fluctuate around these deterministic equivalents as Gaussian random variables, and we will characterize their asymptotic 
mean and variance.


\section{Asymptotic behavior of $\hat{d}_{M}$} \label{sec:firstorder}

The objective of this section is to analyze the asymptotic behavior of the statistic $\hat{d}_{M}(1,2)$ in its most general form. To simplify the notation we will drop the argument $(1,2)$ of the statistic, always keeping in mind that it is built from the two SCMs $\hat{\mathbf{R}}_1$, $\hat{\mathbf{R}}_2$. We will make the following assumptions:

\noindent \textbf{(As1)} For $j\in\{1,2\}$ and $k=1,\ldots,N_{j}$ the observations
$\mathbf{y}_{j}(k)$ are all independent and identically distributed and can be expressed as
\[
\mathbf{y}_{j}(k)=\mathbf{R}_{j}^{1/2}\mathbf{x}_{j}(k)
\]
where $\mathbf{R}_{j}$ is an Hermitian positive definite matrix and $\mathbf{x}_{j}(k)$ is a vector of i.i.d. random entries with zero mean and unit variance. We will consider a boolean variable $\varsigma$ that will indicate whether the observations are real or complex valued. If $\varsigma=1$ the observations are real valued, whereas $\varsigma=0$ indicates they are complex circularly symmetric. It is worth noting that, in this work, we assume no structural model over the covariance matrices $\mathbf{R}_j$. Similarly, the asymptotic descriptors presented in the remainder of this section do not rely on any special data distribution (e.g., Gaussian assumption).

\noindent \textbf{(As2)} The different eigenvalues of $\mathbf{R}_{j}$ are denoted $0<\gamma_{1}^{(j)}<\ldots<\gamma_{\bar{M}_{j}}^{(j)}$ ($j \in \{1,2\}$)\ and have multiplicity $K_{1}^{(j)},\ldots,K_{\bar{M}_{j}}^{(j)}$, where $\bar{M}_{j}$ is the total number of distinct eigenvalues. All these quantities may vary with $M$ \ but we always have $\inf_{M}\gamma_{1}^{(j)}>0$ and $\sup_{M}
\gamma_{\bar{M}_j}^{(j)}<\infty$.

\noindent \textbf{(As3)} The quantities $N_{1}$ and $N_{2}$ depend on $M$, that is
$N_{1}=N_{1}(M)$ and $N_{2}=N_{2}(M)$. Furthermore, when $M\rightarrow\infty$ we have, for $j\in\{1,2\}$, $N_{j}(M)\rightarrow\infty$ in a way that $M/N_{j}
\rightarrow c_{j}$ for some constant $0<c_{j}<\infty$ such that $c_{j}\neq1$.

Assumption \textbf{(As1)} is quite standard, except for the fact that we introduce the Boolean variable $\varsigma$. 
Assumption \textbf{(As2)} points out that the eigenvalues of the true covariance matrices $\mathbf{R}_j$, $j \in \{1,2\}$ may behave freely as the dimensions of the matrix grows, as long as they fluctuate in a compact interval of the positive real axis independent of $M$.
Finally, it is worth pointing out that \textbf{(As3)} explicitly excludes the case $c_j = 1$, mainly because addressing this case requires more elaborate technical tools that are out of the scope of this paper. 

In order to analyze the behavior of $\hat{d}_{M}$ under the above assumptions,
we will build upon random matrix theory tools. 
To begin with, let us consider the function of complex
variable $\omega_{j}\left(  z\right)$, given by one of the solutions to the
polynomial equation
\begin{equation}
z=\omega_{j}\left(  z\right)  \left(  1-\frac{1}{N_{j}}\sum_{m=1}^{\bar{M}_{j}}K_{m}^{(j)}\frac{\gamma_{m}^{(j)}}{\gamma_{m}^{(j)}-\omega_{j}\left(
z\right)  }\right)  . \label{eq:defw(z)}
\end{equation}
More specifically, if $z\in\mathbb{C}^{+}$ (upper complex semiplane),
$\omega_{j}\left(  z\right)  $ is the only solution of the above equation located in $\mathbb{C}
^{+}$. If $z\in\mathbb{C}^{-}$(lower complex semiplane), $\omega_{j}\left(
z\right)  $ is the only solution in $\mathbb{C}^{-}$. Finally, if $z$ is real
valued, $\omega_{j}\left(  z\right)  $ is defined as the only real valued
solution such that
\begin{equation}
\frac{1}{N_{j}}\sum_{m=1}^{\bar{M}_{j}}K_{m}^{(j)}\left(  \frac{\gamma_{m}^{(j)}}{\gamma_{m}^{(j)}-\omega_{j}\left(  z\right)  }\right)  ^{2}<1.
\label{eq:conditionwrealvalued}
\end{equation}
The following result is well known in the literature (it can easily be obtained from \cite{Silverstein95}). It basically states that we can asymptotically replace the random matrices $f_j^{(l)}(\hat{\mathbf{R}}_j)$ in (\ref{eq:definitionhatd}) with some deterministic matrices (usually referred to as deterministic equivalents) in the sense that the difference between the two functionals will converge to zero with probability one. 

\begin{theorem}
\label{th:asymptoticDetEq}Let $\mathbf{A}_{M}$ denote a sequence of
deterministic $M\times M$ matrices with bounded spectral norm. For $z\in\mathbb{C}^{+}
$, consider the resolvents
\begin{align*}
\mathbf{\hat{Q}}_{j}(z)  &  =\left(  \mathbf{\hat{R}}_{j}-z\mathbf{I}
_{M}\right)  ^{-1} \text{and} & 
\mathbf{Q}_{j}(\omega)  =\left(  \mathbf{R}_{j}-\omega\mathbf{I}
_{M}\right)  ^{-1}
\end{align*}
for $j\in\{1,2\}$. Under $\textbf{(As1)}-\textbf{(As3)}$ we have
\[
\frac{1}{M}\mathrm{tr}\left[  \mathbf{A}_{M}\mathbf{\hat{Q}}_{j}(z)\right]
-\frac{\omega_{j}\left(  z\right)  }{z}\frac{1}{M}\mathrm{tr}\left[
\mathbf{A}_{M}\mathbf{Q}_{j}\left(  \omega_{j}\left(  z\right)  \right)
\right]  \rightarrow0
\]
almost surely. 

\end{theorem}
Theorem~\ref{th:asymptoticDetEq} states that the random resolvent $\hat{\mathbf{Q}}_j(z)$ 
asymptotically behaves as deterministic matrix, given by 
\begin{equation}
\mathbf{\bar{Q}}_{j}\left(  z\right)
=\frac{\omega_{j}\left(  z\right)  }{z}\mathbf{Q}_{j}\left(  \omega_{j}\left(
z\right)  \right)  . \label{eq:defQbar}
\end{equation}
By connecting the functions  $f_j^{(l)}(\hat{\mathbf{R}}_j)$ with the resolvent $\hat{\mathbf{Q}}_j(z)$ (see \textbf{(As4)} below) we will be able to formulate an equivalent property for it.

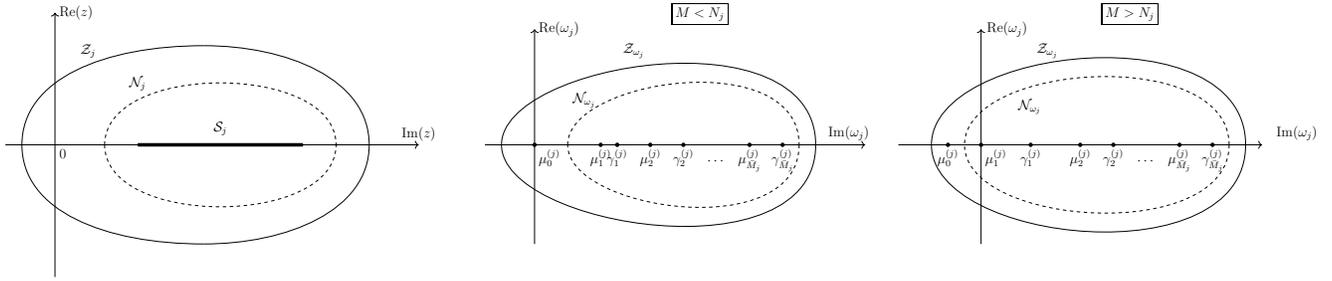
\begin{figure*}[t]
    \centering
    \resizebox{2\columnwidth}{!}{

\tikzset{
    ultra thin/.style= {line width=0.1pt},
    very thin/.style=  {line width=0.2pt},
    thin/.style=       {line width=0.4pt},
    semithick/.style=  {line width=0.6pt},
    thick/.style=      {line width=0.8pt},
    very thick/.style= {line width=1.2pt},
    ultra thick/.style={line width=3pt}
}

\begin{tikzpicture}[
thick,
font=\large,]

\draw [->] (-4,1.5) -- (8.5,1.5);
\draw [->] (-2.5,-2.5) -- (-2.5,5.5);
\draw [dashed] (-1,1.5) .. controls (-1,4) and (6,4) .. (6,1.5);
\draw [dashed] (-1,1.5) .. controls (-1,-1) and (6,-1) .. (6,1.5);
\draw [style = ultra thick] (0,1.5) -- (5,1.5);
\node at (2.5,2) {$\mathcal{S}_j$};
\node[below right] at (-2.5,1.5) {$0$};
\draw (7,1.5) .. controls (7,3.5) and (4.5,4.5) .. (2,4.5);
\draw (7,1.5) .. controls (7,-0.5) and (4.5,-1.5) .. (2,-1.5);
\draw (2,4.5) .. controls (-1,4.5) and (-3.5,3.5) .. (-3.5,1.5);
\draw (2,-1.5) .. controls (-1,-1.5) and (-3.5,-0.5) .. (-3.5,1.5);
\node [right] at (-2.5,5.5) {$\mathrm{Re}(z)$};
\node  [above] at (8.5,1.5) {$\mathrm{Im}(z)$};
\node [above] at (0,3) {$\cnz{j}$};
\node [above] at (-1.5,4) {$\cz{j}$};

\node [above] at (17,5) {$\boxed{M<N_j}$};
\node [above] at (30,5) {$\boxed{M>N_j}$};

\draw [->] (12,-1.5) -- (12,5);
\draw [->]  (10.5,1.5) -- (22,1.5);
\draw [fill=black] (12,1.5) ellipse (0.05 and 0.05) node [below right] {$\mu^{(j)}_0$};
\draw [fill=black] (14,1.5) ellipse (0.05 and 0.05) node [below] {$\mu^{(j)}_{1}$};
\draw [fill=black] (14.5,1.5) ellipse (0.05 and 0.05) node [below] {$\gamma^{(j)}_{1}$};
\draw [fill=black] (15.5,1.5) ellipse (0.05 and 0.05) node [below] {$\mu^{(j)}_{2}$};
\draw [fill=black] (16.5,1.5) ellipse (0.05 and 0.05) node [below] {$\gamma^{(j)}_{2}$};
\draw [fill=black] (18.5,1.5) ellipse (0.05 and 0.05) node [below] {$\mu^{(j)}_{\bar{M}_j}$};
\draw [fill=black] (19.5,1.5) ellipse (0.05 and 0.05) node [below] {$\gamma^{(j)}_{\bar{M}_j}$};
\node at (17.5,1) {$\cdots$};
\draw (20.5,1.5) .. controls (20.5,5.5) and (11,4) .. (11,1.5);
\draw (20.5,1.5) .. controls (20.5,-2.5) and (11,-1) .. (11,1.5);
\draw [dashed] (13,1.5) .. controls (13,3.5) and (20,4.5) .. (20,1.5);
\draw [dashed] (20,1.5) .. controls (20,-1.5) and (13,-0.5) .. (13,1.5);
\node [above] at (13.5,2.5) {$\cnz{\omega_j}$};
\node [above] at (15,4) {$\cz{\omega_j}$};
\node [right] at (12,5) {$\mathrm{Re}(\omega_j)$};
\node  [above] at (21.5,1.5) {$\mathrm{Im}(\omega_j)$};

\draw [->] (25.5,-1.5) -- (25.5,5);
\draw [->] (23,1.5) -- (34,1.5);

\draw [fill=black] (25.5,1.5) ellipse (0.05 and 0.05) node [below right] {$\mu^{(j)}_1$};
\draw [fill=black] (24.5,1.5) ellipse (0.05 and 0.05) node [below] {$\mu^{(j)}_{0}$};
\draw [fill=black] (27,1.5) ellipse (0.05 and 0.05) node [below] {$\gamma^{(j)}_{1}$};
\draw [fill=black] (28.5,1.5) ellipse (0.05 and 0.05) node [below] {$\mu^{(j)}_{2}$};
\draw [fill=black] (29.5,1.5) ellipse (0.05 and 0.05) node [below] {$\gamma^{(j)}_{2}$};
\draw [fill=black] (31.5,1.5) ellipse (0.05 and 0.05) node [below] {$\mu^{(j)}_{\bar{M}_j}$};
\draw [fill=black] (32.5,1.5) ellipse (0.05 and 0.05) node [below] {$\gamma^{(j)}_{\bar{M}_j}$};
\node [below right] at (26.5,3) {$\cnz{\omega_j}$};
\node [above] at (27.5,4) {$\cz{\omega_j}$};
\node [right] at (25.55,5) {$\mathrm{Re}(\omega_j)$};
\node  [above] at (35.05,1.5) {$\mathrm{Im}(\omega_j)$};
\node at (30.5,1) {$\cdots$};
\draw [dashed] (33,1.5) .. controls (33,4.5) and (25,4) .. (25,1.5);
\draw [dashed] (33,1.5) .. controls (33,-1.5) and (25,-1) .. (25,1.5);
\draw (33.5,1.5) .. controls (33.5,5.5) and (24,4.5) .. (24,1.5);
\draw (33.5,1.5) .. controls (33.5,-2.5) and (24,-1.5) .. (24,1.5);
\end{tikzpicture}
    }
    \caption{Illustration of the two possible choices of contour considered in this paper (left), and their transformation with the mapping $z \mapsto \omega_j(z)$, denoted as $\cnz{\omega_j} = \omega_j(\cnz{j})$ and $\cz{\omega_j} = \omega_j(\cz{j})$. In the center, the configuration obtained when $N_j>M$ (oversampled case). On the right, the configuration when $N_j<M$ (undersampled case).}
    \label{fig:contours}
\end{figure*}

More specifically, let us see how this result can be immediately used to analyze the asymptotic behavior of
$\hat{d}_{M}$ in (\ref{eq:definitionhatd}) under some additional assumptions on the quantities \smash{$f_j^{(l)}(\hat{\mathbf{R}}_j)$}. 
Indeed, consider the interval $\mathcal{S}_{j}=[\theta_{j}^{-},\theta_{j}
^{+}]$, where
\begin{equation}
\theta_{j}^{-}=\inf_{M}\gamma_{1}^{(j)}\left(  1-\sqrt{c_{M,j}}\right)
^{2},\,\theta_{j}^{+}=\sup_{M}\gamma_{\bar{M}_j}^{(j)}\left(  1+\sqrt
{c_{M,j}}\right)  ^{2} \label{eq:limitsOuterSupport}
\end{equation}
where $c_{M,j}={M}/{N_{j}}$.
Now, according to \cite{Bai98}, all the positive eigenvalues of $\mathbf{\hat{R}}_{j}$ belong to $\mathcal{S}_{j}$ with probability one for all $M$ sufficiently large. Using this property, we will assume that the functions $f_{j}^{(l)}(\mathbf{\hat{R}}_{j})  $ accept the following representation:

\noindent \textbf{(As4)} For $j\in\{1,2\}$ and $l=1,\ldots,L,$ the quantity $f_{j}
^{(l)} (\mathbf{\hat{R}}_{j})  $ can be expressed as
\begin{equation}
f_{j}^{(l)}\left(  \mathbf{\hat{R}}_{j}\right)  =\frac{1}{2\pi\mathrm{j}}
\oint\nolimits_{C_{j}^{(l)}}f_{j}^{(l)}(z)\mathbf{\hat{Q}}_{j}(z)dz
\label{eq:fRhatasInt}
\end{equation}
with probability one for all large $M$, where $C_{j}^{(l)}$ is a negatively oriented simple closed contour enclosing $\mathcal{S}_{j}$ and not crossing zero and where, with some abuse of notation, \smash{$f_{j}^{(l)}(z)$} denotes a complex function analytic on an open set including $C_{j}^{(l)}$. 
In particular, when \smash{$f_{j}^{(l)}(z)$} is analytic on the set $\mathcal{S}_j$ we can see the left hand side of (\ref{eq:fRhatasInt}) as the result of the Cauchy integration formula. 
However, the result that we present below is more general and also applies to the general situation in which $f_{j}^{(l)}(z)$ may present singularities inside $\mathcal{S}_{j}$.
\begin{remark}
    All the contours in this paper should be understood to be negatively (counter-clockwise) oriented. Furthermore, we will denote as $\cz{j}$ a general contour that encloses $\mathcal{S}_j \cup \{0\}$ and as $\cnz{j}$ a contour that encloses $\mathcal{S}_j$ but not $\{0\}$. See Figure \ref{fig:contours} (left) for an illustration of these contours.  
\end{remark}

The following result follows from a direct application of Theorem
\ref{th:asymptoticDetEq} together with the dominated convergence theorem.

\begin{proposition}
\label{prop:generalEstimator}Under \textbf{(As1)}-\textbf{(As4)} and for any given sequence of deterministic matrices $\mathbf{A}_M$ with bounded norm,
\begin{equation}
\frac{1}{M} \mathrm{tr} \left[ f_{j}^{(l)}(  \mathbf{\hat{R}}_{j}) \mathbf{A}_M \right]  - \frac{1}{2\pi\mathrm{j}
}\oint\nolimits_{C_{j}^{(l)}} \frac{f_{j}^{(l)}(z)}{M} \mathrm{tr} \left[ \mathbf{\bar{Q}}_{j}(z) \mathbf{A}_M \right] dz \rightarrow 0
\label{eq:asymptEqF}
\end{equation}
where $\mathbf{\bar{Q}}_{j}(z)$ is as defined in (\ref{eq:defQbar}).
\end{proposition}
\begin{IEEEproof}
Consider $M$ sufficiently large and the probability set for which (\ref{eq:fRhatasInt})
holds true (a set which, by $\textbf{(As4)}$, has probability one). We can
write
\begin{multline*}
f_{j}^{(l)}\left(  \mathbf{\hat{R}}_{j}\right)  -\frac{1}{2\pi\mathrm{j}}
\oint\nolimits_{C_{j}^{(l)}}f_{j}^{(l)}(z)\mathbf{\bar{Q}}_{j}(z)dz = \\
=\frac
{1}{2\pi\mathrm{j}}\oint\nolimits_{C_{j}^{(l)}}f_{j}^{(l)}(z)\left[
\mathbf{\hat{Q}}_{j}(z)-\mathbf{\bar{Q}}_{j}(z)\right]  dz.
\end{multline*}
Moreover, by omitting the dependence on $M$ in $\mathbf{A}_M$, we obtain
\begin{multline*}
\left\vert \oint\nolimits_{C_{j}^{(l)}}f_{j}
^{(l)}(z)\frac{1}{M}\mathrm{tr}\left[  \mathbf{A}\left(  \mathbf{\hat{Q}}
_{j}(z)-\mathbf{\bar{Q}}_{j}(z)\right)  \right]  dz\right\vert \leq \\ \leq\sup_{z\in
C_{j}^{(l)}}\left\vert f_{j}^{(l)}(z)\right\vert \oint\nolimits_{C_{j}^{(l)}
}\left\vert \frac{1}{M}\mathrm{tr}\left[  \mathbf{A}\left(  \mathbf{\hat{Q}
}_{j}(z)-\mathbf{\bar{Q}}_{j}(z)\right)  \right]  \right\vert \left\vert
dz\right\vert
\end{multline*}
where, obviously, \smash{$\sup_{z\in C_{j}^{(l)}} | f_{j}^{(l)}(z) | 
<\infty$} because of its analyicity. We know from Theorem
\ref{th:asymptoticDetEq} that $M^{-1}\mathrm{tr}[\mathbf{A}(\mathbf{\hat{Q}
}_{j}(z)-\mathbf{\bar{Q}}_{j}(z))]\rightarrow0$ almost surely for all fixed
$z\in C_{j}\cap\mathbb{C}^{+}$. However, $\mathbf{\hat{Q}}_{j}(z)$ and
$\mathbf{\bar{Q}}_{j}(z)$ are analytic functions on an open subset including
$C_{j}^{(l)}$, and from  Lemma \ref{lemma:UnifBounds} in Appendix
\ref{sec:lemmasbounds} we have that $\sup_{M}\sup_{z\in C_{j}^{(l)}} \Vert \mathbf{\hat{Q}
}_{j}(z) \Vert  < \infty $ almost surely and $\sup_{M}\sup_{z\in C_{j}^{(l)}}\left\Vert
\mathbf{\bar{Q}}_{j}(z)\right\Vert <\infty$. Hence, it follows that $\sup_{z\in
C_{j}^{(l)}}M^{-1}\mathrm{tr}[\mathbf{A}(\mathbf{\hat{Q}}_{j}(z)-\mathbf{\bar{Q}}
_{j}(z))]\rightarrow0$ as a direct application of Montel's theorem.

\end{IEEEproof}
Proposition \ref{prop:generalEstimator} above shows that $f_{j}^{(l)} (\mathbf{\hat{R}}_{j})$ has a deterministic asymptotic equivalent but
fails to provide a closed form expression for this quantity. In practice, one needs to particularize the integral in (\ref{eq:asymptEqF}) to the specific choices of $f_{j}^{(l)}(z)$ in order to obtain a closed form expression for the corresponding asymptotic equivalent. In order to do that, it turns out to be particularly useful to use the change of variable proposed in \cite{Mestre08tsp} using the invertible map $z\mapsto\omega_{j}(z)$. 
Let us denote by
$
z_{j}\left(  \omega\right)  =\omega ( 1-\frac{1}{N_{j}}
\mathrm{tr}[\mathbf{R}_j\mathbf{Q}(\omega)]
)
$ the inverse of this map 
and let $z_{j}^{\prime}\left(  \omega\right)$ denote its derivative, that is
\begin{equation}
z_{j}^{\prime}\left(  \omega\right)  =1-\Gamma_{j}\left(  \omega\right)
\label{eq:derivativez(w)}
\end{equation}
where we have defined $\Gamma_{j}(  \omega)  ={N_{j}^{-1}}\mathrm{tr}[ \mathbf{R}_j^2\mathbf{Q}^2 _j(\omega)]$.
All this notation allows us to write, by direct application of a change of variables,
\begin{multline}
\frac{1}{2\pi\mathrm{j}}\oint\nolimits_{C_{j}^{(l)}}f_{j}^{(l)}(z)\mathbf{\bar
{Q}}_{j}(z)dz= \\ =\frac{1}{2\pi\mathrm{j}}\oint\nolimits_{C_{\omega_{j}}^{(l)}}
f_{j}^{(l)}(z_{j}\left(  \omega\right)  )\frac{\omega}{z_{j}\left(
\omega\right)  }\mathbf{Q}_{j}(\omega)z_{j}^{\prime}\left(  \omega\right)
d\omega\label{eq:integralChangeVars}
\end{multline}
where $C_{\omega_{j}}^{(l)}=\omega_{j}(C_{j}^{(l)}) $. Now, we can directly use Proposition~\ref{prop:generalEstimator}
to establish that $\hat{d}_{M}$ in (\ref{eq:definitionhatd}) has a
deterministic asymptotic equivalent.

\begin{corollary}
Under \textbf{(As1)}-\textbf{(As4)} we have $\hat{d}_{M}-\bar{d}_{M} \rightarrow 0$ with
probability one, where
\begin{multline}
\bar{d}_{M}=\sum_{l=1}^{L}\frac{1}{\left(  2\pi\mathrm{j}\right)  ^{2}}\oint\nolimits_{C_{1}
^{(l)}}\oint\nolimits_{C_{2}^{(l)}}f_{1}^{(l)}(z_{1})f_{2}
^{(l)}(z_{2}) \times \\ \times \frac{1}{M}\mathrm{tr}\left[  \mathbf{\bar{Q}}_{1}(z_{1}
)\mathbf{\bar{Q}}_{2}(z_{2})\right]  dz_{1}dz_{2}. \label{eq:integral2resolvents}
\end{multline}
\end{corollary}

\begin{IEEEproof}
This will be a direct consequence of Proposition \ref{prop:generalEstimator}
if we are able to establish that $\sup_{M}\Vert f_{j}^{(l)}(
\mathbf{\hat{R}}_{j})  \Vert <\infty$ almost surely for
$j\in\{1,2\}$ and $l=1,\ldots,L$. Now, it follows from the expression in
(\ref{eq:fRhatasInt}) that we have the bound
\[
\left\Vert f_{j}^{(l)}\left(  \mathbf{\hat{R}}_{j}\right)  \right\Vert
\leq\frac{1}{2\pi}\oint\nolimits_{C_{j}}\left\vert f_{j}^{(l)}(z)\right\vert
\left\Vert \mathbf{\hat{Q}}_{j}(z)\right\Vert \left\vert dz\right\vert .
\]
Again $\sup_{z\in C_{j}}\vert f_{j}^{(l)}(z)\vert <\infty$ and
$\sup_{M,z\in C_{j}}\left\Vert \mathbf{\hat{Q}}_{j}(z)\right\Vert
<\infty$ a.s. from Lemma \ref{lemma:UnifBounds} in Appendix~\ref{sec:lemmasbounds}, concluding the proof.
\end{IEEEproof}

In the following subsections we will illustrate how to
solve the double integral in (\ref{eq:integral2resolvents}) in some 
specific cases of interest. The key point will always be the use of
the change of variable stemming from the invertible map $z_j\mapsto
\omega_{j}(z_j)$. As we will see below, this change of variable will
allow us to obtain a closed form expression for the asymptotic 
equivalent $\bar{d}_{M}$. We will mainly focus on the Euclidean distance
(which is given for illustrative purposes), the symmetrized KL divergence
and the subspace distance. A similar derivation could be done for the
log-Euclidean distance using the integration techniques developed in 
\cite{pereira_icassp23}, although details are omitted due to space limitations. 

\subsection{Euclidean distance} \label{sec:Euclideandist1stOrder}

The euclidean distance in (\ref{eq:definitionhatd}) can be expressed according to \textbf{(As4)} with all the contours enclosing $\{0\}$ and
\[
\sum_{l=1}^{L}f_{1}^{(l)}(z_{1})f_{2}^{(l)}(z_{2})=\left(  z_{1}-z_{2}\right)
^{2}=z_{1}^{2}-2z_{1}z_{2}+z_{2}^{2}.
\]
Consequently, to evaluate $\bar{d}_{M}^{E}$ we need to solve the integral in
(\ref{eq:integralChangeVars}) with $f_{j}^{(l)}(z)=z$ and $f_{j}
^{(l)}(z)=z^{2}$, where $C_j^{(l)}=\cz{j}$ for all $l$. Beginning with
$f_{j}^{(l)}(z)=z$ we see that, using the change of variables above, we can
write
\[
\frac{1}{2\pi\mathrm{j}}\oint\nolimits_{\cz{j}}z\mathbf{\bar{Q}
}_{j}(z)dz=\frac{1}{2\pi\mathrm{j}}\oint\nolimits_{\cz{\omega_j}}\omega\mathbf{Q}_{j}(\omega)z_{j}^{\prime}\left(  \omega\right)  d\omega
\]
where $z_{j}^{\prime}\left(  \omega\right)  $ is as in
(\ref{eq:derivativez(w)}) and $\cz{\omega_j} = \omega_j (\cz{j})$. Now the right hand side of the above integrand only has singularities at the eigenvalues $\gamma_{m}^{(j)}$, $m=1,\ldots,\bar{M}_{j}$, which are all enclosed by $\cz{\omega_{j}}$ \cite{Mestre08tsp}, see further Figure \ref{fig:contours}. We can
therefore enlarge the contour $\cz{\omega_{j}}$ and apply a second change of variables $\zeta=\omega^{-1}$ in a way that $\zeta(\cz{\omega_{j}})$ encloses zero and no other singularity.
The corresponding integral takes the form
\[
\frac{1}{2\pi\mathrm{j}}\oint\nolimits_{\cz{\omega_{j}}}z\mathbf{\bar{Q}}_{j}(z)dz=\frac{-1}{2\pi\mathrm{j}}
\oint\nolimits_{\zeta(\cz{\omega_{j}}) }
\frac{z_{j}^{\prime} (\zeta^{-1})\left(  \zeta\mathbf{R}_{j}-\mathbf{I}_{M}\right)^{-1} }{\zeta^{2}} d\zeta.
\]
The only singularity of the integrand corresponds to a second order pole at $\{0\}$, so that observing that $\cz{\omega_{j}}$ is positively oriented and computing the residue at $\zeta=0$ we have
\[
\frac{1}{2\pi\mathrm{j}}\oint\nolimits_{\cz{\omega_{j}}}z\mathbf{\bar{Q}
}_{j}(z)dz=\mathbf{R}_{j}.
\]
Following exactly the same approach one finds that
\[
\frac{1}{2\pi\mathrm{j}}\oint\nolimits_{\cz{\omega_{j}}}z^{2}
\mathbf{\bar{Q}}_{j}(z)dz=\mathbf{R}_{j}^{2}+\left(  \frac{1}{N_{j}
}\mathrm{tr}\mathbf{R}_{j}\right)  \mathbf{R}_{j}
\]
and consequently
\[
\bar{d}_{M}^{E}=\frac{1}{M}\mathrm{tr}\left[  \left(  \mathbf{R}
_{1}-\mathbf{R}_{2}\right)  ^{2}\right]  +\frac{1}{MN_{1}}\mathrm{tr}
^{2}\mathbf{R}_{1}+\frac{1}{MN_{2}}\mathrm{tr}^{2}\mathbf{R}_{2}.
\]
Observe that, as expected, the asymptotic equivalent is different from  the Euclidean distance between the true covariance matrices. This illustrates the fact that distances between SCMs  are generally inconsistent estimators of the corresponding distance  between the true covariance matrices. In some cases, such as in the case of the Euclidean distance, it is possible to modify $\hat{d}_M$ so that it  converges to the Euclidean distance between the true covariances~\cite{pereira_icassp23}.  However, this is not always possible for all distances in the form of~(\ref{eq:definitionhatd}), particularly in the undersampled regime.

\subsection{Symmetrized Kullback-Leibler divergence}

The symmetrized KL distance in (\ref{eq:defJefferiesdivergence}) and its generalization to the undersampled regime (using pseudo-inverses) 
can both be expressed
as in \textbf{(As4)} with
\[
\sum_{l=1}^{L}f_{1}^{(l)}(z_{1})f_{2}^{(l)}(z_{2})=\frac{1}{2}\frac{z_{2}
}{z_{1}}+\frac{1}{2}\frac{z_{1}}{z_{2}}-1
\]
where in the first two terms the contours do not enclose zero, whereas they do
in the last term. Therefore, to find the asymptotic equivalent $\bar{d}
_{M}^{KL}$ we need to evaluate the integral in (\ref{eq:integralChangeVars}) with
$f_{j}^{(l)}(z)=z$  and $f_{j}^{(l)}(z)=z^{-1}$, in both
cases assuming that the corresponding contour does not contain zero (so that $C_j^{(l)}=\cnz{j}$ for all $l$). As the first case has already been considered before, let us
study the integral for $f_{j}^{(l)}(z)=z^{-1}$ and observe that we can
particularize (\ref{eq:integralChangeVars}) to
\begin{equation*}
\frac{1}{2\pi\mathrm{j}}
\oint\nolimits_{\cnz{j}}
\frac{\mathbf{\bar{Q}}_{j}(z)dz}{z}
=
\frac{1}{2\pi\mathrm{j}}
\oint\nolimits_{\cnz{\omega_{j}}}\frac{ \mathbf{Q}_{j}(\omega) (1-\Gamma_{j}\left(  \omega\right) ) }{\omega\left(
1-\frac{\mathrm{tr}[\mathbf{R}_j\mathbf{Q}_j(\omega)]}{N_{j}}  \right)  ^{2}}d\omega
\end{equation*}
where $\cnz{\omega_{j}}=\omega_{j}(\cnz{j})$. Now, let $\mu_{0}^{(j)}<\ldots<\mu_{\bar{M}_{j}}^{(j)}$ denote the solutions to the equation
\begin{equation}
\mu\left(  1-\frac{1}{N_{j}}\mathrm{tr}\left[\mathbf{R}_j\mathbf{Q}_j(\mu)\right]\right)  =0.
\label{eq:def_mu}
\end{equation}
It can be easily  shown \cite{Mestre08tsp} that $\mu_{0}^{(j)}<0$ when $N_{j}<M$ (undersampled case) and $\mu_{0}^{(j)}=0$ when $N_{j}>M$ (oversampled case), see further Figure \ref{fig:contours}. Either way, it turns out that $\mu_{0}^{(j)}$ is the only root in the above set that is not enclosed by
$\cnz{\omega_{j}}$, and therefore all the singularities of the above integrand fall inside $\cnz{\omega_{j}}$ except for a potential singularity at $\mu_{0}^{(j)}$.

Using the above, we can enlarge $\cnz{\omega_{j}}$ in the above
integral so that it encloses $\mu_{0}^{(j)}$ if we then add the corresponding
residue, which turns out to be equal to
$
{\mathbf{R}_{j}\mathbf{Q}_{j}^{2}(\mu_{0}^{(j)})}(1-\Gamma_{j}(
\mu_{0}^{(j)}))^{-1}
$
(notice that $\cnz{\omega_{j}}$ is negatively oriented). Once this has been
evaluated, we can apply the change of variables $\zeta=\omega_{j}^{-1}$, leading
to
\begin{multline*}
\frac{1}{2\pi\mathrm{j}}\oint\nolimits_{\cnz{j}}z^{-1}
\mathbf{\bar{Q}}_{j}(z)dz=\frac{\mathbf{R}_{j}\mathbf{Q}_{j}^{2}(\mu_{0}
^{(j)})}{1-\Gamma_{j}\left(  \mu_{0}^{(j)}\right)  } +\\ +\frac{1}{2\pi\mathrm{j}
}\oint\nolimits_{C_{0}}\frac{1-\Gamma_{j}\left(  \zeta^{-1}\right)
}{\left(  1-\frac{1}{N_{j}}\sum_{m=1}^{\bar{M}_{j}}K_{m}^{(j)}\frac{\gamma
_{m}^{(j)}\zeta}{\gamma_{m}^{(j)}\zeta-1}\right)  ^{2}}\zeta^{-1}
\mathbf{Q}_{j}(\zeta^{-1})d\zeta
\end{multline*}
where now $C_{0}$ is a negatively oriented contour enclosing zero and no other singularity. Noting that the second term of the above expression is zero (the integrand presents an removable singularity at zero), we can conclude that
\[
\frac{1}{2\pi\mathrm{j}}\oint\nolimits_{\cnz{j}}z^{-1}
\mathbf{\bar{Q}}_{j}(z)dz=\frac{\mathbf{R}_{j}\mathbf{Q}_{j}^{2}(\mu_{0}
^{(j)})}{1-\Gamma_{j}\left(  \mu_{0}^{(j)}\right)  }.
\]
With all the above intermediate results, it follows directly that
\[
\bar{d}_{M}^{KL}=\frac{
\mathrm{tr}\left[  \mathbf{R}_{1}
\mathbf{Q}_{1}^{2}(\mu_{0}^{(1)})\mathbf{R}_{2}\right]  }{2M\left(1-\Gamma_{1}\left(\mu_{0}^{(1)}\right) \right)  }
+
\frac{\mathrm{tr}\left[  \mathbf{R}
_{2}\mathbf{Q}_{2}^{2}(\mu_{0}^{(2)})\mathbf{R}_{1}\right]  }{2M\left(1-\Gamma
_{2}\left(  \mu_{0}^{(2)}\right) \right) }-1.
\]
It is particularly interesting to note that in the oversampled situation, that
is when $N_{1},N_{2}>M$, we have $\mu_{0}^{(1)}=\mu_{0}^{(2)}=0$ and
therefore
\[
\bar{d}_{M}^{KL}=\frac{1}{2M}\left(  \frac{N_{1}\mathrm{tr}\left[
\mathbf{R}_{1}^{-1}\mathbf{R}_{2}\right]  }{N_{1}-M}+\frac{N_{2}
\mathrm{tr}\left[  \mathbf{R}_{2}^{-1}\mathbf{R}_{1}\right]  }{N_{2}
-M}\right)  -1.
\]
\subsection{Subspace distance}
The subspace distance also responds to the form in (\ref{eq:definitionhatd})-(\ref{eq:asymptEqF}) with all the functions $f_{j}^{(l)}(z)=1$ and none of
the contours enclose $\{0\}$, that is $C_j^{(l)} = \cnz{j}$ for all $l$, so that
\[
\hat{d}_{M}^{SS}=\frac{N_{1}+N_2}{M}+\frac{1}{2\pi^2}\oint\nolimits_{\cnz{1}}\oint
\nolimits_{\cnz{2}}\frac{\mathrm{tr}[  \mathbf{\hat{Q}%
}_{1}(z_{1})\mathbf{\hat{Q}}_{2}(z_{2})] }{M} dz_{1}dz_{2}.
\]
Using the above integration technique we directly see that
\begin{align*}
&\frac{1}{2\pi\mathrm{j}}\oint\nolimits_{\cnz{j}}\mathbf{\bar{Q}%
}_{j}(z)dz=\frac{1}{2\pi\mathrm{j}}\oint\nolimits_{\cnz{\omega_{j}}}\omega\mathbf{Q}_{j}(\omega)\frac{1-\Gamma_{j}\left(  \omega\right)
}{z_{j}\left(  \omega\right)  }d\omega
\\
&=\mu_{0}^{(j)}\mathbf{Q}_{j}(\mu_{0}^{(j)})+\frac{1}{2\pi\mathrm{j}}%
\oint\nolimits_{\cz{\omega_j}}\omega\mathbf{Q}_{j}(\omega
)\frac{1-\Gamma_{j}\left(  \omega\right)  }{z_{j}\left(  \omega\right)
}d\omega
\\
&=\mathbf{R}_{j}\mathbf{Q}_{j}(\mu_{0}^{(j)})
\end{align*}
and consequently
\[
\bar{d}_{M}^{SS}=\frac{N_{1}}{M}+\frac{N_{2}}{M}-\frac{2}{M}\mathrm{tr}\left[
\mathbf{R}_{1}\mathbf{Q}_{1}(\mu_{0}^{(1)})\mathbf{R}_{2}\mathbf{Q}_{2}%
(\mu_{0}^{(2)})\right]  .
\]%
Obviously, this distance only makes sense in the undersampled regime (otherwise, one cannot possibly define the original subspaces).

\section{Asymptotic fluctuations} \label{sec:secondorder}

The results in the previous two sections were derived for functionals of two SCMs, namely $\hat{\mathbf{R}}_j$ for $j \in \{1,2\}$. In this section, however, we would like
to characterize the \emph{joint} asymptotic statistical behavior of multiple distances between pairs in
a collection of SCMs. To that effect, let $\mathcal{J}$ denote a discrete
set that indexes the set of available SCMs, namely $\hat{\mathbf{R}}_j$ for $j \in \mathcal{J}$. 
Consider a collection of $R$ functionals between pairs of SCMs, where
$R\in \mathbb{N}$ is a fixed integer. We will gather these quantities into an $R$-dimensional
column vector 
\[
\hat{\mathbf{d}}_M = [\hat{d}_M(i_1,j_1),\ldots, \hat{d}_M(i_R,j_R)]^T
\]
where $i_r,j_r \in \mathcal{J}$, for $r=1,\ldots,R$. By construction, some of these indices 
may coincide (so that a SCM may appear in several of these terms). 
We consider that the 
assumptions in the previous sections, which were originally formulated for $\mathcal{J} = \{1,2\}$
are trivially extended for the more general index set $\mathcal{J}$. Obviously, it follows
from previous sections that $\hat{\mathbf{d}}_M - \bar{\mathbf{d}}_M \rightarrow 0$, where 
$\bar{\mathbf{d}}_M$ contains the asymptotic deterministic equivalents $\bar{d}_M(i_r,j_r)$ 
for $r=1,\ldots,R$, which are trivially defined by extending the index set in (\ref{eq:integral2resolvents}). 

Our objective here is to analyze the fluctuations of $\hat{\mathbf{d}}_M$ around its asymptotic equivalent $\bar{\mathbf{d}}_M$. More specifically, let us
consider the following normalized random vector $\hat{\boldsymbol{\zeta}}_{M}=M\left(  \hat{\mathbf{d}}_{M}-\bar{\mathbf{d}}_{M}\right)$. The $r$th entry of this vector can be expressed as 
\begin{align}
\hat{\zeta}_M & (r) =\sum_{l=1}^{L} \frac{1}{4\pi^2}\oint\nolimits_{C_{i_r}^{(l)}}\oint\nolimits_{C_{j_r}^{(l)}
}f_{i_r}(z_{i_r}) f_{j_r}(z_{j_r}) \times \label{eq:definitionzetaM} \\ & \times \mathrm{tr}\left[  \mathbf{\hat{Q}}_{i_r}
(z_{i_r})\mathbf{\hat{Q}}_{j_r}(z_{j_r})-\mathbf{\bar{Q}}_{i_r}(z_{i_r})\mathbf{\bar{Q}
}_{j_r}(z_{j_r})\right]  dz_{i_r}dz_{j_r}. \nonumber
\end{align}
We will establish a central limit theorem (CLT) on $\hat{\boldsymbol{\zeta}}_{M}$ that will basically state that it asymptotically behaves as a Gaussian random vector with a certain mean and variance that are provided next. 

In order to introduce the relevant quantities that will describe the asymptotic mean and variance, we need some notation. For a given $M\times M$ deterministic matrix $\mathbf{A}$, we denote
\begin{equation}
    \Omega_{j}\left(  \omega;\mathbf{A}\right)  =\mathbf{A}+\phi_{j}\left(
\omega;\mathbf{A}\right)  \mathbf{I}_{M}\label{eq:definitionOmega}
\end{equation}
where $\phi_{j}\left(  \omega;\mathbf{A}\right) $ is the scalar function
\begin{equation}
\phi_{j}\left(  \omega;\mathbf{A}\right)  =\frac{\omega}{1-\Gamma_{j}(\omega
)}\frac{1}{N_{j}}\mathrm{tr}\left[  \mathbf{R}_{j}\mathbf{Q}_{j}%
^{2}\left(  \omega\right)  \mathbf{A}\right]  .\label{eq:definitionphij}
\end{equation}

We define the asymptotic (second order) mean of $\hat{\boldsymbol{\zeta}}_{M}$ as an $R$-dimensional
column vector $\boldsymbol{\mathfrak{m}}_{M} = [\mathfrak{m}_{M} (1), \ldots, \mathfrak{m}_{M} (R)]^T$, where
\begin{multline}
\mathfrak{m}_{M} (r)  = \sum_{l=1}^{L} \frac{-\varsigma}{ 4\pi^2} \oint\nolimits_{C_{i_r}^{(l)}}\oint\nolimits_{C_{j_r}^{(l)}}\frac{\omega_{i_r}}{z_{i_r}
}\frac{\omega_{j_r}}{z_{j_r}} f^{(l)}_{i_r}(z_{i_r})f^{(l)}_{j_r}(z_{j_r}) \times \\
\times \mathfrak{m}_{i_r,j_r}\left(  \omega_{i_r},\omega_{j_r}\right)  dz_{i_r}dz_{j_r}
\label{eq:asymptMean}
\end{multline}
where we have introduced the bivariate function
$\mathfrak{m}_{i,j}\left(  \omega_{i},\omega_{j}\right)     =\mathfrak{m}_{i}\left(
\omega_{i},\mathbf{Q}_{j}\left(  \omega_{j}\right)  \right)  +\mathfrak{m}
_{j}\left(  \omega_{j},\mathbf{Q}_{i}\left(  \omega_{i}\right)  \right)$
with 
\begin{equation}
    \mathfrak{m}_{j}\left(  \omega_{j},\mathbf{A}\right)     =\frac{1}{N_{j}%
}\frac{\mathrm{tr}\left[  \mathbf{\mathbf{R}}_{j}^{2}\mathbf{Q}_{j}^{3}\left(
\omega_{j}\right)  \Omega_{j}\left(  \omega_{j};\mathbf{A}\right)  \right]
}{1-\Gamma_{j}(\omega_{j})} \label{eq:defmofomegaA}
\end{equation} 
(valid for $i,j \in \mathcal{J}$)
and where we have used the shorthand notation $\omega_{j}=\omega_{j}(z_{j})$.

The asymptotic variance matrix of the normalized vector $\hat{\boldsymbol{\zeta}}_M$ is denoted as the $R \times R$ matrix $\boldsymbol{\Sigma}_M$. For $r,s =1,\ldots,R$, the $(r,s)$th entry of this matrix quantifies the asymptotic cross-covariance between the $r$th normalized functional $M(\hat{d}_M(i_r,j_r)-\bar{d}_M(i_r,j_r))$ and the $s$th normalized functional $M(\hat{d}_M(i_s,j_s)-\bar{d}_M(i_s,j_s))$. By definition of these functionals and noting that the different SCMs are statistically independent, we see that the $(r,s)$th entry of the $\boldsymbol{\Sigma}_M$, denoted by $\boldsymbol{\Sigma}_M(r,s)$, will be zero unless one of the involved SCM appears both in the $r$th and the $s$th functionals ($\hat{d}_M(i_r,j_r)$ and $\hat{d}_M(i_s,j_s)$). 
Now, observe that we always have $i_r \neq j_r$ (i.e. $\hat{\mathbf{R}}_{i_r} \neq \hat{\mathbf{R}}_{j_r}$) and $i_s \neq j_s$ (i.e. $\hat{\mathbf{R}}_{i_s} \neq \hat{\mathbf{R}}_{j_s}$ ), because each functional is built from two \emph{different} SCMs. Therefore, $\boldsymbol{\Sigma}_M (r,s)$ will only be different from zero if at least one of these conditions hold: $i_r=i_s$, $j_r=i_s$, $i_r = j_s$ or $j_r = j_s$, corresponding to the cases where 1\textsuperscript {st}/2\textsuperscript{nd}/1\textsuperscript{st}/2\textsuperscript{nd} SCM in the $r$th functional is equal to 1\textsuperscript{st}/1\textsuperscript{st}/2\textsuperscript{nd}/2\textsuperscript{nd} SCM in the $s$th functional respectively. 
When at least one of these conditions is true, we have\footnote{The integrals are over the contours $C_{i_r}^{(l_r)} \times C_{j_r}^{(l_r)} \times C_{i_s}^{(l_s)}  \times C_{j_s}^{(l_s)}$. Also, note that we are writing again $\omega_{j}=\omega_{j}(z_{j})$ and $\omega_{j}^{\prime}=\omega_{j}(z_{j}^{\prime})$.}
\begin{align}
\boldsymbol{\Sigma}_M(r,s) &= \sum_{l_r,l_s=1}^{L} \frac{1+\varsigma}{\left(  2\pi\mathrm{j}\right)  ^{4}} 
{\oint \oint \oint \oint}  \frac{\omega_{i_r}\omega_{j_r}}{z_{i_{r}} z_{j_{r}}}\frac{\omega
_{i_s}^{\prime}\omega_{j_s}^{\prime}}{z_{i_s}^{\prime}z_{j_s}^{\prime}} \times \label{eq:asymptVariance}  \\
& \times f_{i_r}^{(l_r)}(z_{i_r}) f_{j_r}^{(l_r)}(z_{j_r}) f_{i_s}^{(l_s)}(z_{i_s}^{\prime}) f_{j_s}^{(l_s)}(z_{j_s}^{\prime}) \times \nonumber \\  
& \times
\sigma_{i_r,j_r,i_s,j_s}
^{2}\left(  \omega_{i_r},\omega_{j_r},\omega_{i_s}^{\prime},\omega_{j_s}^{\prime
}\right)  dz_{i_r}dz_{j_r} dz_{i_s}^{\prime}dz_{j_s}^{\prime} \nonumber 
\end{align}
where we have defined $\sigma^{2}_{i,j,m,n} ( \omega_{i},\omega_{j},\omega_{m}^{\prime},\omega_{n}
^{\prime})$ as
\begin{multline}
\sigma^{2}_{i,j,m,n}\left(  \omega_{i},\omega_{j},\omega_{m}^{\prime},\omega_{n}
^{\prime}\right) =\\ 
= \sigma_{i}^{2}\left(  \omega_{i},\omega_{m}^{\prime};\mathbf{Q}_{j}\left(
\omega_{j}\right)  ,\mathbf{Q}_{n}(\omega_{n}^{\prime})\right) \delta_{i=m}  + \\ 
+ \sigma_{j}^{2}\left(  \omega_{j},\omega_{n}^{\prime};\mathbf{Q}_{i}\left(
\omega_{i}\right)  ,\mathbf{Q}_{m}(\omega_{m}^{\prime})\right) \delta_{j=n} +\\
+ \sigma_{i}^{2}\left(  \omega_{i},\omega_{n}^{\prime};\mathbf{Q}_{j}\left(
\omega_{j}\right)  ,\mathbf{Q}_{m}(\omega_{m}^{\prime})\right) \delta_{i=n} + \\
+ \sigma_{j}^{2}\left(  \omega_{j},\omega_{m}^{\prime};\mathbf{Q}_{i}\left(
\omega_{i}\right)  ,\mathbf{Q}_{n}(\omega_{n}^{\prime})\right) \delta_{j=m} +\\
+ \varrho_{i,j}(\omega_i,\omega_j,\omega_m^{\prime},\omega_n^{\prime}) \delta_{i=m}\delta_{j=n} +\\
+ \varrho_{i,j}(\omega_i,\omega_j,\omega_n^{\prime},\omega_m^{\prime}) \delta_{i=n}\delta_{j=m} \label{eq:Sigma2}
\end{multline}
with the following additional definitions. The functions $\sigma_{j}^{2}\left(  \omega,\omega^{\prime};\mathbf{A},\mathbf{B}\right) $ for $j \in \mathcal{J} $ are defined in (\ref{eq:omegabivariate}) at the top of the  page, where we have introduced the quantities 
\[
\Gamma_{j}(\omega,\omega^{\prime})=\frac{1}{N_{j}}\mathrm{tr}\left[
\mathbf{R}_{j}^{2} \mathbf{Q}_j\left(\omega\right) \mathbf{Q}_j\left(\omega^{\prime}\right)
\right].
\]
\begin{figure*}[t]
    \centering
    \begin{align}
\sigma_{j}^{2}\left(  \omega,\omega^{\prime};\mathbf{A},\mathbf{B}\right)   &
=\frac{1}{1-\Gamma_{j}(\omega,\omega^{\prime})}\frac{1}{N_{j}}\mathrm{tr}
\left[  \mathbf{R}_{j}\mathbf{Q}_{j}\left(  \omega\right)  \mathbf{Q}
_{j}\left(  \omega^{\prime}\right)  \Omega_{j}\left(  \omega;\mathbf{A}
\right)  \mathbf{R}_{j}\mathbf{Q}_{j}\left(  \omega\right)  \mathbf{Q}
_{j}\left(  \omega^{\prime}\right)  \Omega_{j}\left(  \omega^{\prime
};\mathbf{B}\right)  \right] \nonumber \\
&  +\frac{1}{\left(  1-\Gamma_{j}(\omega,\omega^{\prime})\right)  ^{2}}
\frac{1}{N_{j}}\mathrm{tr}\left[  \mathbf{R}_{j}^{2}\mathbf{Q}_{j}^{2}\left(
\omega\right)  \mathbf{Q}_{j}\left(  \omega^{\prime}\right)  \Omega_{j}\left(
\omega;\mathbf{A}\right)  \right]  \frac{1}{N_{j}}\mathrm{tr}\left[
\mathbf{R}_{j}^{2}\mathbf{Q}_{j}\left(  \omega\right)  \mathbf{Q}_{j}
^{2}\left(  \omega^{\prime}\right)  \Omega_{j}\left(  \omega^{\prime
};\mathbf{B}\right)  \right] \label{eq:omegabivariate}
\end{align} 
     \hrulefill
\end{figure*}
Regarding $\varrho_{i,j}(\omega_i,\omega_j,\omega_i^{\prime},\omega_j^{\prime}) $, they are defined as 
\begin{multline}
\varrho_{i,j}(\omega_i,\omega_j,\omega_i^{\prime},\omega_j^{\prime}) = \\
=
\frac{\mathrm{tr}^{2}\left[  \mathbf{R}_{i}\mathbf{Q}_{i}\left(  \omega
_{i}\right)  \mathbf{Q}_{i}\left(  \omega_{i}^{\prime}\right)  \mathbf{R}
_{j}\mathbf{Q}_{j}\left(  \omega_{j}\right)  \mathbf{Q}_{j}\left(  \omega
_{j}^{\prime}\right)  \right]  }{N_{i}N_{j}\left(  1-\Gamma_{i}(\omega
_{i},\omega_{i}^{\prime})\right)  \left(  1-\Gamma_{j}(\omega_{j},\omega
_{j}^{\prime})\right)  }.  \label{eq:defvarrho}
\end{multline}
We have now all then necessary notation to introduce the main result of this section, which is a CLT on the random vector $\hat{\boldsymbol{\zeta}}_M$. In order to establish the result, we assume that the observations are Gaussian distributed. A more general CLT could be derived in more general statistical settings, but we prefer to focus on the Gaussian case for simplicity.

\begin{theorem}
\label{theorem:CLT}
Assume that \textbf{(As1)}-\textbf{(As4)} hold for the extended set of SCMs and that the observations are
Gaussian distributed. If the minimum eigenvalue of the matrix $\boldsymbol{\Sigma}_M$ is bounded away from zero, that is $$\liminf_{M\rightarrow\infty}\lambda_{\min}(\boldsymbol{\Sigma}_M)>0$$ then $\boldsymbol{\Sigma}_M^{-1/2} (\hat{\boldsymbol{\zeta}}_{M}-\boldsymbol{\mathfrak{m}}_{M})$ converges in law to a multivariate standard Gaussian distribution. 
\end{theorem}

\begin{IEEEproof}
See supplementary material.
\end{IEEEproof}
The above result can be used to approximate the behavior of the functionals of SCMs $\hat{\mathbf{d}}_M$ for finite values of $M, \{N_j\}_{j\in \mathcal{J}}$. Indeed, one can asymptotically approximate $\hat{\mathbf{d}}_M$ as a multivariate Gaussian random variable with mean value $\bar{\mathbf{d}}_M + \boldsymbol{\mathfrak{m}}_M/M$ and covariance matrix $\boldsymbol{\Sigma}_M/{M^2}$. This will be of fundamental help in order to study the performance of these functionals in specific problems, such as the covariance clustering problem (see Section \ref{sec:numEvaluation} for details). 

In what follows we will see how this general result particularizes to some of the distances between SCM that have been introduced before. For illustrative purposes and due to limited space, we only provide a detailed proof for the $\hat{d}_M^\mathrm{E}$ (in Appendix B). Also, to simplify the exposition, we assume $R=1$ and $\mathcal{J} = \{1,2\}$, the results being trivially extrapolated to the more general case of $R>1$ and multiple SCM. The CLT for the other distances are obtained following the same integration techniques (final expressions are therefore given without proof). Alternatively, the evaluation of the integrals in (\ref{eq:integral2resolvents}), (\ref{eq:asymptMean}) and (\ref{eq:asymptVariance}) can be carried out numerically.

\subsection{Squared Euclidean distance}
By direct evaluation of the integrals in (\ref{eq:asymptMean})-(\ref{eq:asymptVariance}) when $g(z_1,z_2) = (z_1-z_2)^2$ one can establish that $\mathfrak{m}_M$ takes the form 
\begin{equation}
\mathfrak{m}_{M}=\varsigma\left(\frac{1}{N_{1}}\mathrm{tr}\left[
\mathbf{\mathbf{R}}_{1}^{2}\right]  +\frac{1}{N_{2}}\mathrm{tr}\left[
\mathbf{\mathbf{R}}_{2}^{2}\right]\right)
    \label{eq:secondOrderMeanEuclidean}
\end{equation}
whereas $\Sigma_M$ particularizes  to
\begin{align}
    \frac{\Sigma_{M}}{1+\varsigma}  &=
 2\left(  \frac{\mathrm{tr}[\mathbf{R}_1^2]}{N_{1}}  \right)^2
 + 
 2\left(  \frac{\mathrm{tr}[\mathbf{R}_2^2]}{N_{2}}  \right)^2
 +\frac{4\mathrm{tr}^2[\mathbf{R}_1\mathbf{R}_2] }{N_{1}N_{2}} 
 \nonumber
 \\
 &+ \frac{4}{N_1}\mathrm{tr}\left[\left(\mathbf{R}_1\Psi_1\right)^2\right] + \frac{4}{N_2}\mathrm{tr}\left[\left(\mathbf{R}_2\Psi_2\right)^2\right]
 \label{eq:varianceEuclidean}
\end{align}  
where we have introduced $\Psi_j = (\mathbf{R}_1 - \mathbf{R}_2) + (M/N_j)\mathbf{I}_M$, for $j \in \{1,2\}$. All the terms are obviously positive (uniformly in $M$), so that one can easily see that the variance is uniformly bounded away from zero. This shows that the CLT holds for the Euclidean distance between SCMs.

\subsection{Symmetrized KL metric}

In this case, we need to consider the integrals in
(\ref{eq:asymptMean})-(\ref{eq:asymptVariance}) for $g(z_{1},z_{2}
)=z_{1}/(2z_{2})+z_{2}/(2z_{1})-1$ where only the contours of the last
integral enclose $\{0\}$. One can show that the second order mean takes the form
\begin{equation}
    \mathfrak{m}_{M}=\varsigma \sum_{\substack{i, j\in\{1,2\} \\ i \neq j}} \frac{\left. d [\omega_i \mathfrak{m}_i(\omega_i,\mathbf{R}_j)] / d\omega_i \right\vert_ {\omega_i = \mu_0^{(i)} }} {2\left(  1-\Gamma_{i}(  \mu_{0}^{(i)})  \right)  } 
    \label{eq:secondOrderMeanKL}
\end{equation}
where we recall that $\mathfrak{m}_j(\omega,\mathbf{A})$ is defined in (\ref{eq:defmofomegaA}).

Regarding the variance, one can proceed in a similar way to show that
\begin{multline}
\frac{\Sigma_{M}}{1+\varsigma}   =\frac{\left.  {\partial^{2}\left[  \omega
_{1}\omega_{1}^{\prime}\Upsilon_{11}\left(  \omega_{1},\omega_{1}^{\prime
}\right)  \right]  } /{\partial\omega_{1}\partial\omega_{1}^{\prime}}\right\vert
_{\omega_{1}=\omega_{1}^{\prime}=\mu_{0}^{(1)}}}{4\left(  1-\Gamma_{1}(  \mu
_{0}^{(1)})  \right)  ^{2}}\\
  +\frac{\left.  {\partial^{2}\left[  \omega_{2}\omega_{2}^{\prime
}\Upsilon_{22}\left(  \omega_{2},\omega_{2}^{\prime}\right)  \right]
}/{\partial\omega_{2}\partial\omega_{2}^{\prime}}\right\vert  _{\omega_{2}
=\omega_{2}^{\prime}=\mu_{0}^{(2)}}}{4\left(  1-\Gamma_{2}(  \mu_{0}^{(2)})
\right)  ^{2}}\\
  +\frac{\left.
{\partial^{2}\left[  \omega_{1}\omega_{2}\Upsilon_{12}\left(  \omega
_{1},\omega_{2}\right)  \right]  }/{\partial\omega_{1}\partial\omega_{2}%
}\right\vert _{\omega_{1}=\mu_{0}^{(1)},\omega_{2}=\mu_{0}^{(2)}}}{2\left(  1-\Gamma_{1}( \mu_{0}^{(1)})
\right)  \left(  1-\Gamma_{2}(  \mu_{0}^{(2)})  \right)  }
\label{eq:varianceKL}
\end{multline}
where we have defined
\begin{multline*}
\Upsilon_{11}\left(  \omega_{1},\omega_{1}^{\prime}\right)     = \frac{\mathrm{tr}^{2}\left[  \mathbf{R}
_{2}\mathbf{R}_{1}\mathbf{Q}_{1}\left(  \omega_{1}\right)  \mathbf{Q}
_{1}\left(  \omega_{1}^{\prime}\right)  \right]  }{N_{1}N_{2}\left(
1-\Gamma_{1}(\omega_{1},\omega_{1}^{\prime})\right) }  +  \\+\sigma
_{1}^{2}\left(  \omega_{1},\omega_{1}^{\prime};\mathbf{R}_{2},\mathbf{R}
_{2}\right) + \frac{1}{N_{2}}\mathrm{tr}\left[  \mathbf{R}_{2}\mathbf{Q}_{1}\left(  \omega_{1}\right)  \mathbf{R}_{2}\mathbf{Q}_{1}\left(  \omega_{1}^{\prime}\right)  \right] 
\end{multline*}
where $\Upsilon_{22}\left(  \omega_{2},\omega_{2}^{\prime}\right)$ is defined equivalently but interchanging the two indexes ($1 \leftrightarrow 2$) and where
\begin{multline*}
 \Upsilon_{12}\left(  \omega_{1},\omega_{2}\right)    = 
 \frac{1}{N_{1}N_{2}}\mathrm{tr}^{2}\left[  \mathbf{R}_{1}\mathbf{Q}_{1}\left(  \omega_{1}\right)  \mathbf{R}_{2}\mathbf{Q}_{2}\left(  \omega_{2}\right)  \right] 
 \\
 -\frac{1}{N_{1}}\mathrm{tr}\left[  \mathbf{R}_{1}\mathbf{Q}_{1}\left(  \omega_{1}\right)\mathbf{Q}_{2}\left(  \omega_{2}\right) \mathbf{R}_{1}\mathbf{Q}_{1}\left(
\omega_{1}\right)  \Omega_{1}\left(  \omega_{1};\mathbf{R}_{2}\right)  \right]
\\
  -\frac{1}{N_{2}}\mathrm{tr}\left[  \mathbf{R}_{2}\mathbf{Q}_{2}\left(
\omega_{2}\right)  \mathbf{Q}_{1}\left(  \omega_{1}\right)  \mathbf{R}%
_{2}\mathbf{Q}_{2}\left(  \omega_{2}\right)  \Omega_{2}\left(  \omega
_{2};\mathbf{R}_{1}\right)  \right] .
\end{multline*}

The expression of the second order mean and variance can be significantly simplified in the oversampled case, where we will always have $\mu_0^{(1)}=\mu_0^{(2)}=0$. In this situation, the second order mean
particularizes to 
\[
\mathfrak{m}_M  = \frac{\varsigma}{2}\left[\frac{N_1 \mathrm{tr}[\mathbf{R}_2\mathbf{R}_1^{-1}]}{(N_1-M)^2}+
\frac{N_2 \mathrm{tr}[\mathbf{R}_1 \mathbf{R}_2^{-1}]}{(N_2-M)^2}\right]
\]
whereas the variance takes the simple form
\[
\frac{\Sigma_M^2}{1+\varsigma} = \frac{N_1^2 \Upsilon_{11}(0,0)}{4(N_1-M)^2}  + \frac{N_2^2 \Upsilon_{22}(0,0)}{4(N_2-M)^2}
+ \frac{N_1 N_2 \Upsilon_{12}(0,0)}{2(N_1-M)(N_2-M)}
\]
where 
\[
\Upsilon_{11}(0,0) = \frac{N_1+N_2-M}{N_2(N_1-M)}\left[ \mathrm{tr}[(\mathbf{R}_1^{-1}\mathbf{R}_2)^2] + \frac{\mathrm{tr}^2[\mathbf{R}_1^{-1}\mathbf{R}_2]}{N_1-M} \right] 
\]
with $\Upsilon_{22}(0,0)$ equivalently defined by swapping indexes ($1 \leftrightarrow 2$), and where 
\[
\Upsilon_{12}(0,0)= \frac{M^2}{N_1 N_2}- \frac{M}{N_1}- \frac{M}{N_2}.
\]
The fact that $\liminf_{M}\sigma_M^2 >0 $ is easy to see in the
oversampled case, by simply using the fact that $\alpha\mathrm{tr}\mathbf{A} + \beta\mathrm{tr}\mathbf{A}^{-1}>2M\sqrt{\alpha\beta}$ for $\alpha,\beta>0$ and positive definite $\mathbf{A}$. 
Indeed, by using this inequality with $\mathbf{A} = (\mathbf{R}_1^{-1}\mathbf{R}_2)^2$ and noting that the terms of the form $\mathrm{tr}^2(\cdot)$ are positive, we see that 
\begin{multline*}
    \frac{\Sigma_M}{1+\varsigma} > \frac{N_1 N_2}{2(N_1-M)(N_2-M)}
\times \\
 \times \left( 
\frac{M(N_1+ N_2- M)}{\sqrt{(N_1-M)(N_2-M) N_1 N_2}}
+\frac{M^2}{N_1 N_2} 
- \frac{M}{N_1}- \frac{M}{N_2}\right).
\end{multline*}
Next, we observe that $(N_1+N_2-M)/(N_1-M)>1+N_2/N_1$, so that
\[
\frac{M^2(N_1+N_2-M)^2}{(N_1-M)(N_2-M)N_1N_2}>\left( \frac{M}{N_1}+\frac{M}{N_2} \right)^2
\]
and we can conclude that $\sigma_M^2 > M^2/2/(N_1-M)/(N_2-M)$, which is bounded away from zero.

\subsection{Squared subspace distance}

In the case of the squared subspace distance, 
the second order mean takes the form
\begin{multline}
\mathfrak{m}_{M}   =-2\mu_{0}^{(1)}\mathfrak{m}_{1}(  \mu_{0}
^{(1)},\mathbf{R}_{2}\mathbf{Q}_{2}(  \mu_{0}^{(2)}))  \\
 -2\mu_{0}^{(2)}\mathfrak{m}_{2}(  \mu_{0}^{(2)},\mathbf{R}%
_{1}\mathbf{Q}_{1}(  \mu_{0}^{(1)}) )  .
    \label{eq:2ndordermeanProj}
\end{multline}
whereas the asymptotic variance can be written as
\begin{align}
\frac{\Sigma_{M}}{1+\varsigma} &  =4(  \mu_{0}^{(1)})
^{2}\sigma_{1}^{2}(  \mu_{0}^{(1)},\mu_{0}^{(1)};\mathbf{R}_{2}%
\mathbf{Q}_{2}(  \mu_{0}^{(2)})  ,\mathbf{R}_{2}\mathbf{Q}%
_{2}(  \mu_{0}^{(2)})  )  \nonumber \\
&  +4(  \mu_{0}^{(2)})  ^{2}\sigma_{2}^{2}(  \mu_{0}^{(2)}%
,\mu_{0}^{(2)};\mathbf{R}_{1}\mathbf{Q}_{1}(  \mu_{0}^{(1)})
,\mathbf{R}_{1}\mathbf{Q}_{1}(  \mu_{0}^{(1)})  )  \nonumber \\
&  +4(  \mu_{0}^{(1)}\mu_{0}^{(2)})  ^{2}\frac{\mathrm{tr}%
^{2}[  \mathbf{R}_{1}\mathbf{Q}_{1}^{2}(  \mu_{0}^{(1)})
\mathbf{R}_{2}\mathbf{Q}_{2}^{2}(  \mu_{0}^{(2)}) ]  }%
{N_{1}N_{2}(  1-\Gamma_{1}(\mu_{0}^{(1)}))  (  1-\Gamma
_{2}(\mu_{0}^{(2)}))  }. \label{eq:varianceProj}
\end{align}
Close examination of the expression of the variance reveals that, since $\mathbf{Q}_j(\mu_0^{(j)})$ is positive definite, the first two terms of (\ref{eq:varianceProj}) are non-negative. Moreover, it is easy to see that $|\mu_0^{(j)}| \geq (M/N_j-1)\gamma_{\bar{M}_j}^{(j)}$, which is bounded away from zero. A direct application of Lemma \ref{lemma:UnifBounds} in Appendix~\ref{sec:lemmasbounds} shows that the third term in (\ref{eq:varianceProj}) is bounded away from zero, and hence the CLT holds. 
Finally, we highlight that one can also estimate these quantities by following similar steps as the ones in~\cite{roberto22_icassp}.

\section{Numerical Evaluation} 
\label{sec:numEvaluation}

In order to validate the results presented above, we consider the set of SCMs $\{\hat{\mathbf{R}}_j\}_{j \in \mathcal{J}}$, each of which is obtained over $N_j$ observations and is associated to (possibly distinct) Toeplitz\footnote{Note that the results above are general and do not consider any structure over the covariance matrices. We do so in this numerical evaluation only as a way to generate data in a controlled environment. } covariance matrices $\mathbf{R}_j$ with first rows $[\rho_j^0, \dots, \rho_j^{M-1}]$. 
We will carry out the analysis in two steps.
First, through numerical evaluation, 
we study the validity of the asymptotic statistical description of the family of quantities $\hat{d}_M$ defined in Section~\ref{sec:statistical_model_observations}.
Specifically, we are interested in the particularizations $\hat{d}_M^{E}$, $\hat{d}_M^{KL}$, $\hat{d}_M^{SS}$. 
Second, we employ the derived solutions to assess the quality of the metrics applied to the problem of clustering SCMs. 
We recall that, in the oversampled regime, the 
subspace similarity metric does not exist. Hence,
we deliberately evaluate this metric only in the undersampled regime.

\subsection{Consistency of Asymptotic Descriptors}

\begin{figure}
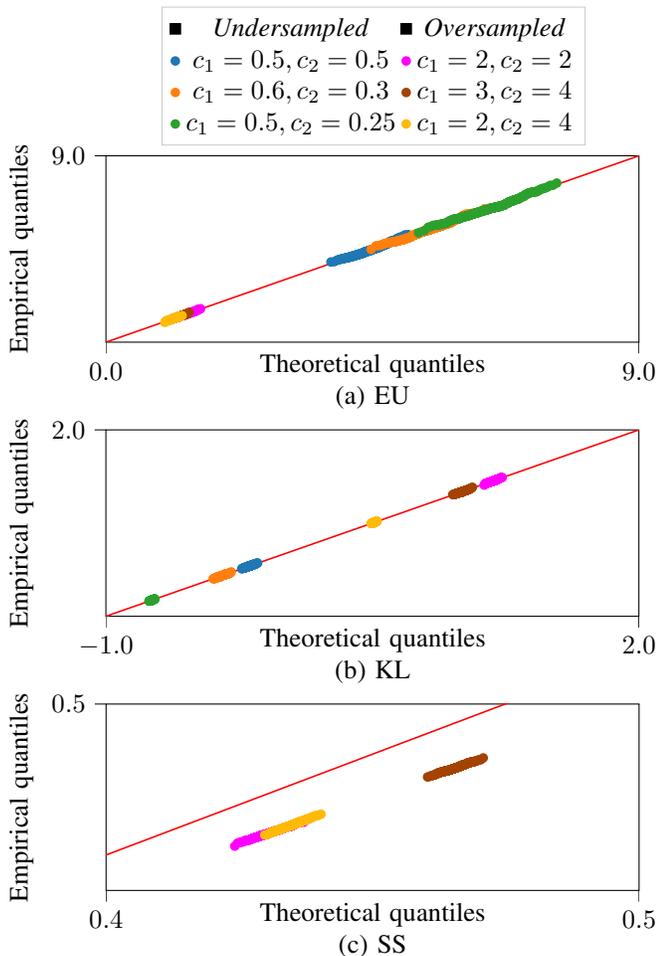

\begin{subfigure}
{\input{figs/qq_plots/many/qq_many_eu}}
 \vspace{-2em}
\end{subfigure}
\hspace{-1em}

\begin{subfigure}
{\input{figs/qq_plots/many/qq_many_kl}}
\vspace{-2em}
\end{subfigure}

\begin{subfigure}
{\input{figs/qq_plots/many/qq_many_pf}}
  \vspace{-2em}
\end{subfigure}

     \caption{
    Quantitle-Quantile plot of empirical distribution and asymptotic descriptors of the different metrics (a) EU, (b) KL and (c) SS.}
\label{fig:results:histograms}
\end{figure}

We start by comparing the asymptotic (theoretical) descriptors from the {theorems} above against the empirical distribution of their respective metrics, namely, Euclidean (EU), Symmetrized Kullback Leibler (KL) and subspace (SS) distances. 
A common graphical method for such a comparisons involves plotting the quantiles of the random variable against its theoretical values. 
Fig.~\ref{fig:results:histograms} illustrates this comparison by presenting multiple quantile-quantile (qq) plots obtained from the comparison of the empirical and theoretical distributions for $\mathcal{J} = \{1,2\}$, $M=200$, $\rho_1=0.7$, $\rho_2=0.8$, and various choices of $c_1, c_2$. Each quantile comparison is represented by a different solid circle in the plot. 
Observe that there is a very good match between the theoretical and the empirical quantile values regardless of the considered metric or of whether $c_1, c_2$ are large or small. This tendency is visually represented by a close alignment between the circles and the solid lines which indicate a perfect match between two distributions. In other words, we have that the random distribution of the distance $\hat{d}_M$ seems to be correctly approximated by 
$$
\hat{d}_M \sim \mathcal{N}\left( \bar{d}_M + \frac{\mathfrak{m}_M}{M}, \frac{\Sigma_M}{M^2} \right).
$$

\subsection{Assessing the Clustering of SCMs}

An interesting consequence of the results above is that they allow us to compare the suitability of distance metrics for a specific clustering problem without the need of collecting any data.  In other words, by using the  asymptotic results above, one is capable of predicting the behavior of a clustering algorithm based solely on the statistics of the underlying functionals.  Consider the problem of clustering  six sample covariance matrices $\hat{\mathbf{R}}_1,\ldots,\hat{\mathbf{R}}_6$, which are pairwise associated to the same covariance matrix. In other words, we have $\rho_1=\rho_2$, $\rho_3=\rho_4$ and $\rho_5=\rho_6$. Furthermore, assume for simplicity that $N_1=N_2$, $N_3=N_4$ and $N_5=N_6$. 
We will define the empirical probability of correct clustering (``Prob. Correct'' in the figures), evaluated over $10^4$ simulations, as the percentage of realizations for which all the distances between SCMs associated to the same covariance matrix (i.e. belonging to the same true cluster) are lower than the minimum distance between SCMs associated to different covariance matrices (i.e. belonging to different true clusters).

\begin{figure}
\centering
\begin{subfigure}
{\hspace{-1.3em}\definecolor{darkgray176}{RGB}{176,176,176}
\definecolor{darkorange25512714}{RGB}{255,127,14}
\definecolor{forestgreen4416044}{RGB}{44,160,44}
\definecolor{lightgray204}{RGB}{204,204,204}
\definecolor{magenta}{RGB}{255,0,255}
\definecolor{orange2551868}{RGB}{255,186,8}
\definecolor{saddlebrown164660}{RGB}{164,66,0}
\definecolor{steelblue31119180}{RGB}{31,119,180}

\begin{tikzpicture} 
\begin{axis}[
hide axis,
width=0.8\linewidth,
height=0.3\linewidth,
at={(0,0)},
xmin=10,
xmax=50,
ymin=0,
ymax=0.4,
legend style={fill opacity=1, draw opacity=1, text opacity=1, draw=lightgray204, legend columns=2,
    anchor=south,
    at={(-4.5, 0)},
    font=\small,
    legend image post style={scale=0.6}
}
]

\addlegendimage{black,mark=square*, opacity=0, only marks}
    \addlegendentry{\quad \textit{Oversampled }};
    
    \addlegendimage{black,mark=square*, opacity=0, only marks}
    \addlegendentry{\quad \textit{}};

\addlegendimage{steelblue31119180,mark=square*, only marks};
\addlegendentry{$c_1=2.0, c_3 = 2.0, c_5 = 2.0$ };

\addlegendimage{darkorange25512714,mark=square*, only marks};
\addlegendentry{$c_1=2.0, c_3 = 3.0, c_5 = 4.0$ };

\addlegendimage{black,mark=square*, opacity=0, only marks}
\addlegendentry{\quad \textit{Undersampled }};

\addlegendimage{black,mark=square*, opacity=0, only marks}
\addlegendentry{\quad \textit{}};

\addlegendimage{forestgreen4416044,mark=square*, only marks};
\addlegendentry{$c_1=0.8, c_3 = 0.8, c_5 = 0.8$ };

\addlegendimage{magenta,mark=square*, only marks};
\addlegendentry{$c_1=0.6, c_3 = 0.6, c_5 = 0.5$ };

\addlegendimage{saddlebrown164660,mark=square*, only marks};
\addlegendentry{$c_1=0.3, c_3 = 0.4, c_5 = 0.7$};

\addlegendimage{orange2551868,mark=square*, only marks};
\addlegendentry{$c_1=0.8, c_3 = 0.3, c_5 = 0.5$};

    \coordinate (legend) at (axis description cs:0.97,0.03);

\end{axis}



\end{tikzpicture}}
\end{subfigure}
\begin{subfigure}
{
\begin{tikzpicture}

\definecolor{darkgray176}{RGB}{176,176,176}
\definecolor{darkorange25512714}{RGB}{255,127,14}
\definecolor{forestgreen4416044}{RGB}{44,160,44}
\definecolor{lightgray204}{RGB}{204,204,204}
\definecolor{magenta}{RGB}{255,0,255}
\definecolor{orange2551868}{RGB}{255,186,8}
\definecolor{saddlebrown164660}{RGB}{164,66,0}
\definecolor{steelblue31119180}{RGB}{31,119,180}

\begin{axis}[
width=0.9\linewidth,
height=0.5\linewidth,
at={(0,0)},
tick align=outside,
tick pos=left,
x grid style={darkgray176},
xlabel={Growing M},
xmin=10, xmax=300,
xtick style={color=black},
y grid style={darkgray176},
ytick={0, 0.5, 1},
ylabel={Prob. Correct},
ymin=0.0, ymax=1.005,
ytick style={color=black},
]
\addplot [semithick, steelblue31119180]
table {%
10 0.4568
12 0.6422
14 0.7967
16 0.8848
18 0.9425
20 0.9739
22 0.9914
24 0.9968
26 0.9985
28 0.9994
30 1
32 1
34 1
36 1
38 1
40 1
42 1
44 1
46 1
48 1
50 1
60 1
70 1
80 1
90 1
110 1
130 1
150 1
250 1
280 1
300 1
};
\addplot [semithick, steelblue31119180, dashed, line width=1.5pt]
table {%
10 0.53714165779038
12 0.71334129744119
14 0.839025515878011
16 0.916698072905281
18 0.959808686742305
20 0.981828880243555
22 0.992274837281217
24 0.996910355613713
26 0.99884680751689
28 0.999602118488036
30 0.999859710671725
32 0.999960265971542
34 0.999985463359595
36 0.999998958863694
38 0.999999603142919
40 0.999998184734672
42 1.0000067188261
44 0.999998917953162
46 1.0000032325901
48 0.999992343217576
50 1.00000334091841
60 1.00000247776635
70 1.00000631393417
80 1.00000081408739
90 1.00000335009169
110 1.00000331805792
130 0.999995567787595
150 0.99999968558888
250 1.00000559843
280 0.999993039154865
300 0.999993042863591
};
\addplot [semithick, darkorange25512714, line width=1.5pt]
table {%
10 0.2627
12 0.3241
14 0.3847
16 0.4506
18 0.5028
20 0.5668
22 0.6298
24 0.6783
26 0.734
28 0.7842
30 0.8182
32 0.8634
34 0.8902
36 0.9104
38 0.9334
40 0.9491
42 0.9613
50 0.9908
60 0.9988
70 0.9997
80 1
90 1
110 1
130 1
150 1
250 1
280 1
300 1
};
\addplot [semithick, darkorange25512714, dashed, line width=1.5pt]
table {%
10 0.236542322857283
12 0.295911928349974
14 0.359662545626917
16 0.426902544987436
18 0.495704592295563
20 0.563982362949388
22 0.629710585533447
24 0.691199085407146
26 0.74712365217305
28 0.796660801592718
30 0.839422977114337
32 0.875452700485092
34 0.905098323053727
36 0.928947250663283
38 0.947723462369679
40 0.96219937099876
42 0.973135809173007
50 0.994237509289145
60 0.999431753897682
70 0.999967060080807
80 0.999998989786473
90 0.999999984032355
110 0.999999999998734
130 1
150 1
250 1
280 1
300 1
};
\addplot [semithick, forestgreen4416044, line width=1.5pt]
table {%
10 0.0134
20 0.0827
30 0.2302
40 0.4136
50 0.6051
60 0.7479
70 0.8469
80 0.916
90 0.9543
110 0.984
130 1
150 1
250 1
280 1
300 1
};
\addplot [semithick, forestgreen4416044, dashed, line width=1.5pt]
table {%
10 0.0181875874712574
20 0.114149473081241
30 0.309136763911336
40 0.531815124369116
50 0.709779209913674
60 0.827786592797324
70 0.900769995469129
80 0.944856228089605
90 0.970730556093544
110 0.993008215567046
130 0.99869544083088
150 0.999810176738592
250 1.00000001340534
280 1.00000000645123
300 1.00000000368574
};
\addplot [semithick, magenta, line width=1.5pt]
table {%
10 0.0118
20 0.0574
30 0.1449
40 0.2458
50 0.3267
60 0.4116
70 0.4972
80 0.5662
90 0.6422
110 0.7557
130 0.844
150 0.908
250 0.994
280 1
300 1
};
\addplot [semithick, magenta, dashed, line width=1.5pt]
table {%
10 0.0101880695732217
20 0.063795494660802
30 0.150054904374431
40 0.241347499519923
50 0.329296050454501
60 0.41759980740716
70 0.504242702430689
80 0.584852294018437
90 0.658507866094966
110 0.779455729386028
130 0.866596919972073
150 0.92435254203721
250 0.998580179456
280 0.999705358606698
300 0.999909002634229
};
\addplot [semithick, saddlebrown164660, line width=1.5pt]
table {%
10 0
12 0
14 0
16 0
18 0
20 0
22 0
24 0
26 0
28 0
30 0
32 0
34 0
36 0
38 0
40 0
42 0
44 0
46 0
48 0
50 0
60 0
70 0
80 0
90 0
110 0
130 0
150 0
250 0
280 0
300 0
};
\addplot [semithick, saddlebrown164660, dashed, line width=1.5pt]
table {%
10 0
12 0
14 0
16 0
18 0
20 0
22 0
24 0
26 0
28 0
30 0
32 0
34 0
36 0
38 0
40 0
42 0
44 0
46 0
48 0
50 0
60 0
70 0
80 0
90 0
110 0
130 0
150 0
250 0
280 0
300 0
};
\addplot [semithick, orange2551868, line width=1.5pt]
table {%
10 0
12 0
14 0
16 0
18 0
20 0
22 0
24 0
26 0
28 0
30 0
32 0
34 0
36 0
38 0
40 0
42 0
44 0
46 0
48 0
50 0
60 0
70 0
80 0
90 0
110 0
130 0
150 0
250 0
280 0
300 0
};
\addplot [semithick, orange2551868, dashed, line width=1.5pt]
table {%
10 0
12 0
14 0
16 0
18 0
20 0
22 0
24 0
26 0
28 0
30 0
32 0
34 0
36 0
38 0
40 0
42 0
44 0
46 0
48 0
50 0
60 0
70 0
80 0
90 0
110 0
130 0
150 0
250 0
280 0
300 0
};
\end{axis}

\end{tikzpicture}}
 \vspace{-1.5em}
\caption*{\hspace{2em}(a) KL}
\end{subfigure}
\hspace{-1em}
\centering
\begin{subfigure}
{
\begin{tikzpicture}

\definecolor{darkgray176}{RGB}{176,176,176}
\definecolor{darkorange25512714}{RGB}{255,127,14}
\definecolor{forestgreen4416044}{RGB}{44,160,44}
\definecolor{lightgray204}{RGB}{204,204,204}
\definecolor{magenta}{RGB}{255,0,255}
\definecolor{orange2551868}{RGB}{255,186,8}
\definecolor{saddlebrown164660}{RGB}{164,66,0}
\definecolor{steelblue31119180}{RGB}{31,119,180}

\begin{axis}[
width=0.9\linewidth,
height=0.5\linewidth,
at={(0,0)},
tick align=outside,
tick pos=left,
x grid style={darkgray176},
xlabel={Growing M},
xmin=10, xmax=300,
xtick style={color=black},
y grid style={darkgray176},
ytick={0, 0.5, 1},
ylabel={Prob. Correct},
ymin=0.0, ymax=1.005,
ytick style={color=black},
]
\addplot [semithick, steelblue31119180, line width=1.5pt]
table {%
10 0.0962
12 0.1682
14 0.2524
16 0.3363
18 0.4367
20 0.503
22 0.5872
24 0.6437
26 0.6928
28 0.7351
30 0.7676
32 0.7981
34 0.8243
36 0.839
38 0.8552
40 0.866
42 0.8815
44 0.8922
46 0.9068
48 0.9137
50 0.9237
60 0.9587
70 0.9769
80 0.9892
90 0.9947
110 0.997
130 1
150 1
250 1
280 1
300 1
};
\addplot [semithick, steelblue31119180, dashed, line width=1.5pt]
table {%
10 0.0571941879840479
12 0.0972887849300375
14 0.14748778538429
16 0.205508122502478
18 0.268806944166815
20 0.334808249344226
22 0.401280966957109
24 0.466440926283683
26 0.528926103850284
28 0.587758369437476
30 0.642331887674685
32 0.692292065883391
34 0.737482343492238
36 0.777888681613201
38 0.813637020512469
40 0.844933769417799
42 0.872021098486224
44 0.895225997146642
46 0.914884934487738
48 0.93141164551403
50 0.945146439027389
60 0.983932284316559
70 0.996042581028891
80 0.999189274699898
90 0.999854711346591
110 0.99999657734914
130 1.0000033279492
150 1.00000273849319
250 1.00000021980956
280 0.999998897526099
300 0.999998677902234
};
\addplot [semithick, darkorange25512714, line width=1.5pt]
table {%
10 0.1502
12 0.2641
14 0.3941
16 0.5351
18 0.6506
20 0.7692
22 0.8358
24 0.8864
26 0.9245
28 0.95
30 0.9604
32 0.9747
34 0.982
36 0.9861
38 0.9917
40 0.9924
42 0.9951
50 0.9991
60 1
70 0.9997
80 1
90 1
110 1
130 1
150 1
250 1
280 1
300 1
};
\addplot [semithick, darkorange25512714, dashed, line width=1.5pt]
table {%
10 0.145373312478936
12 0.242511361369167
14 0.351307908661673
16 0.461774623151615
18 0.566161394873362
20 0.659381539853479
22 0.738841322686657
24 0.803982970553591
26 0.855668343996633
28 0.895452836781744
30 0.925384204950175
32 0.947481270057631
34 0.963479490266108
36 0.974930730436077
38 0.98300117273629
40 0.988567184232611
42 0.992469365208408
50 0.998759363974433
60 0.999923565736544
70 0.999997650526197
80 0.999999931109737
90 0.999999998724995
110 0.999999999999894
130 1
150 1
250 1
280 1
300 1
};
\addplot [semithick, forestgreen4416044, line width=1.5pt]
table {%
10 0.0155999999999999
20 0.0869
30 0.2652
40 0.4528
50 0.6106
60 0.7034
70 0.7683
80 0.8093
90 0.8499
110 0.923
130 0.964
150 0.977
250 1
280 1
300 1
};
\addplot [semithick, forestgreen4416044, dashed, line width=1.5pt]
table {%
10 0.00779664679828173
20 0.0554245738988629
30 0.163659092188847
40 0.316323586213022
50 0.475189480413601
60 0.612203926678161
70 0.719208615098381
80 0.799177964356199
90 0.857921118038776
110 0.932107777602116
130 0.970430128185035
150 0.988557374412555
250 0.999990844230195
280 1.00000212930096
300 1.00000088572988
};
\addplot [semithick, magenta, line width=1.5pt]
table {%
10 0.00619999999999998
20 0.0358000000000001
30 0.1045
40 0.206
50 0.2903
60 0.3674
70 0.4209
80 0.4369
90 0.4617
110 0.4914
130 0.54
150 0.535
250 0.702
280 0.747
300 0.774
};
\addplot [semithick, magenta, dashed, line width=1.5pt]
table {%
10 0.00325272535532848
20 0.0190459542973041
30 0.0546973315553791
40 0.109829451171848
50 0.175415333107338
60 0.239956220480425
70 0.296151447370245
80 0.342093437848444
90 0.379039549147254
110 0.43457057648649
130 0.477169772583932
150 0.514825957952821
250 0.685893756357708
280 0.731656171948087
300 0.760167483482284
};
\addplot [semithick, saddlebrown164660, line width=1.5pt]
table {%
10 0.00209999999999999
20 0.00729999999999997
30 0.021
40 0.0344
50 0.0568
60 0.0806
70 0.102
80 0.1229
90 0.1416
110 0.195
130 0.251
150 0.312
250 0.61
280 0.675
300 0.73
};
\addplot [semithick, saddlebrown164660, dashed, line width=1.5pt]
table {%
10 0.00214423134886154
20 0.013568618054291
30 0.0357282125121259
40 0.0614152685291993
50 0.0846976013568264
60 0.10483461408035
70 0.123481244476291
80 0.142186327994395
90 0.161824629780623
110 0.20504899569085
130 0.25355669951784
150 0.306553409122296
250 0.593470145131282
280 0.670945352517292
300 0.717497113689619
};
\addplot [semithick, orange2551868, line width=1.5pt]
table {%
10 0.00419999999999998
20 0.00870000000000004
30 0.00700000000000001
40 0.003
50 0.00219999999999998
60 0.000900000000000012
70 0.000499999999999945
80 0
90 0
110 0
130 0
150 0
};
\addplot [semithick, orange2551868, dashed, line width=1.5pt]
table {%
10 0.00137222640839976
20 0.0028122882157817
30 0.00295964178262957
40 0.00226461088571918
50 0.00142291964986825
60 0.000757208595507891
70 0.000348931047911872
80 0.000140902906211388
90 5.02950083573805e-05
110 4.55104938544289e-06
130 2.7035657603696e-07
300 0
};
\end{axis}

\end{tikzpicture}}
\vspace{-1.5em}
\caption*{\hspace{2em}(b) EU}
\end{subfigure}
\centering
\begin{subfigure}
{
\begin{tikzpicture}

\definecolor{darkgray176}{RGB}{176,176,176}
\definecolor{green}{RGB}{0,128,0}
\definecolor{lightgray204}{RGB}{204,204,204}


\definecolor{darkorange25512714}{RGB}{255,127,14}
\definecolor{forestgreen4416044}{RGB}{44,160,44}
\definecolor{magenta}{RGB}{255,0,255}
\definecolor{orange2551868}{RGB}{255,186,8}
\definecolor{saddlebrown164660}{RGB}{164,66,0}
\definecolor{steelblue31119180}{RGB}{31,119,180}


\begin{axis}[
width=0.9\linewidth,
height=0.5\linewidth,
at={(0,0)},
tick align=outside,
tick pos=left,
x grid style={darkgray176},
xlabel={Growing M},
xmin=10, xmax=300,
xtick style={color=black},
y grid style={darkgray176},
ytick={0, 0.5, 1},
ylabel={Prob. Correct},
ymin=-0.001, ymax=1.005,
ytick style={color=black}
]
\addplot [semithick, steelblue31119180, line width=1.5pt]
table {%
10 0
300 0
};
\addplot [semithick, steelblue31119180, line width=1.5pt]
table {%
10 0
20 0
30 0
40 0
50 0
60 0
70 0
80 0
90 0
110 0
130 0
150 0
250 0
280 0
300 0
};
\addplot [semithick, darkorange25512714, line width=1.5pt]
table {%
10 0
300 0
};
\addplot [semithick, darkorange25512714, line width=1.5pt]
table {%
10 0
20 0
30 0
40 0
50 0
60 0
70 0
80 0
90 0
110 0
130 0
150 0
250 0
280 0
300 0
};
\addplot [semithick, forestgreen4416044, line width=1.5pt]
table {%
10 0
300 0
};
\addplot [semithick, forestgreen4416044, line width=1.5pt]
table {%
10 0.00207813602122656
20 0.00176768066543326
30 0.00126641049269257
40 0.00076562880305589
50 0.000391499249379737
60 0.000168711931937069
70 6.15998730652831e-05
80 1.90415613081223e-05
90 4.9939065269511e-06
110 2.12288446755737e-07
130 4.85886447253996e-09
150 6.00103043534325e-11
250 2.73154189535511e-24
280 1.64030174080455e-29
300 2.63813743461137e-33
};
\addplot [semithick, magenta, line width=1.5pt]
table {%
10 0
300 0
};
\addplot [semithick, magenta, line width=1.5pt]
table {%
10 0.00289683670699086
20 0.00197601170704333
30 0.000545621115756597
40 7.12652815517513e-05
50 4.85800937595624e-06
60 1.83852931950632e-07
70 4.01074056215342e-09
80 5.10646096549283e-11
90 3.83903420781184e-13
110 4.41375266180181e-18
130 6.21051159894828e-24
150 1.03959452046588e-30
250 1.93135763025071e-78
280 2.74912285874565e-97
300 5.23750403114334e-111
};
\addplot [semithick, saddlebrown164660, line width=1.5pt]
table {%
10 0.0738
20 0.2234
30 0.3757
40 0.5021
50 0.6023
60 0.7
70 0.7652
80 0.8265
90 0.8802
110 0.94
130 0.977
150 0.992
250 1
280 1
300 1
};
\addplot [semithick, saddlebrown164660, dashed, line width=1.5pt]
table {%
10 0.0582210490525985
20 0.210110160353863
30 0.369948780746071
40 0.499673157366745
50 0.602246584941124
60 0.686158514015946
70 0.757156567558164
80 0.817189486987017
90 0.866618142809409
110 0.93564352547497
130 0.972833461860131
150 0.989972231060824
250 0.99998845787191
280 0.9999994310091
300 0.999999924850701
};
\addplot [semithick, orange2551868, line width=1.5pt]
table {%
10 0.1545
20 0.3759
30 0.5107
40 0.6061
50 0.6847
60 0.755
70 0.8124
80 0.869
90 0.902
110 0.959
130 0.983
150 0.994
300 1
};
\addplot [semithick, orange2551868, dashed, line width=1.5pt]
table {%
10 0.141797885083811
20 0.373043652480812
30 0.512651577200283
40 0.605336487375686
50 0.684638853891996
60 0.755405987680144
70 0.816339955851276
80 0.866626980516349
90 0.906394784883728
110 0.958450109277556
130 0.983979460580858
150 0.994638316107211
300 1
};
\end{axis}

\end{tikzpicture}}
    \vspace{-1.5em}
\caption*{\hspace{2em}(c) SS}
\vspace{-0.5em}
\end{subfigure}
     \caption{
    Empirical (solid lines) and theoretical (dashed lines) probability of correct clustering (y-axis) six SCMs into three groups for growing $M$ (x-axis).
    }
\label{fig:results:clustering_SCMs_xaxis_M}
\vspace{-1.5em}
\end{figure}
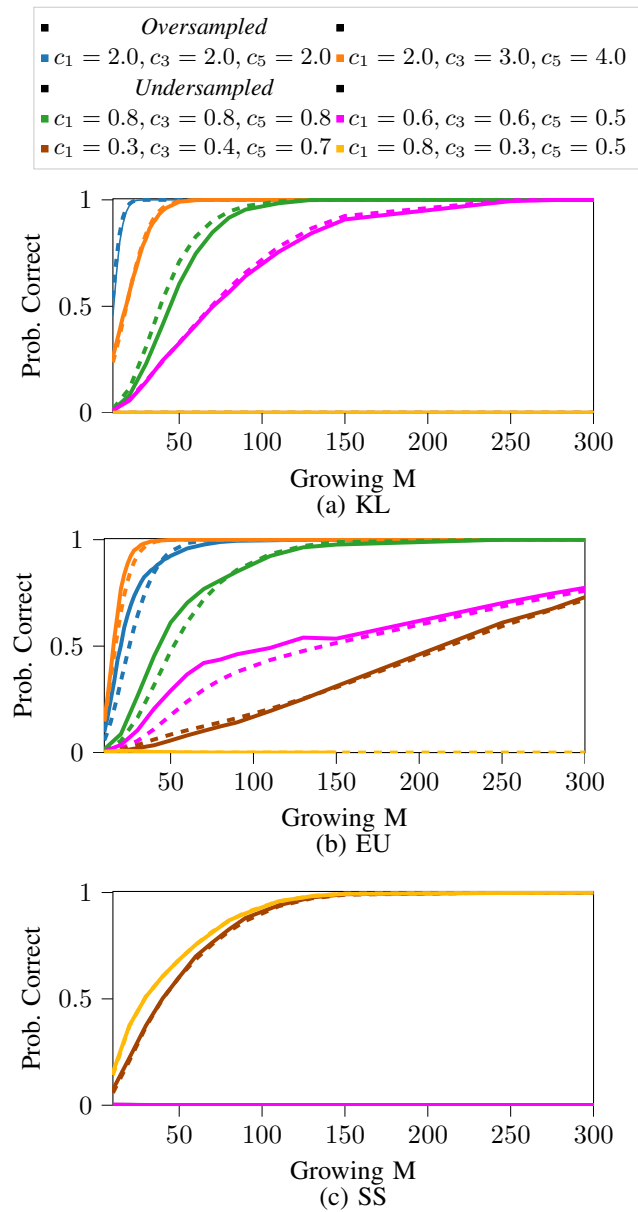

Specifically,  Fig.~\ref{fig:results:clustering_SCMs_xaxis_M} portrays the empirical probability of correct clustering (solid lines), obtained from the $10^4$ different realizations of the SCMs, against the theoretical probability (dashed lines) obtained by directly application of Theorem~\ref{theorem:CLT} for fixed $\rho_1 = \rho_2=0.3, \rho_3=\rho_4 = 0.6, \rho_5 = \rho_6 = 0.9$,  increasing $M, N_j$,  and different values of $c_j, j=1,\ldots, 6$. 
Note that there exists a very good alignment between the empirical and theoretical probabilities of correctly clustering the SCMs, indicating the correctness of our results and a good prediction mechanism. This happens regardless of the scenario (different values of $c_j, j=1, \ldots, 6$) or whether we consider the oversampled or undersampled regime. 

Another interesting analysis comes from the comparing the performance of the different metrics. 
In general terms, none of the metrics (among the ones considered in this work) universally outperforms all the others for every scenario. As an example, in the undersampled regime, the SS distance effectively works when
$c_1=c_2=0.8,c_3=c_4=0.3, c_5=c_6=0.5$ (yellow line) while the symmetrized KL and EU distance,  even as $M$ grows large, fail to accurately cluster the SCMs in this scenario (i.e., have probability zero). In contrast, the SS distance fails 
properly to cluster the SCMs when $c_1=\ldots=c_6=0.8$ or $c_1=\ldots=c_4=0.6, c_5=c_6=0.5$ (green and magenta lines, respectively) while the results for the symmetrized KL and EU distances consistently increase with growing $M$. Naturally, similar behaviors are also present in the oversampled regime. Specifically, in the scenario where $c_1=\ldots=c_6=2.0$ (blue line) the symmetrized KL distance outperforms EU distance, and, conversely, in the case where $c_1=c_2=2.0, c_3=c_4=3.0, c_5=c_6=4.0$ (orange line) the result is reversed, with the EU distance outperforming the symmetrized KL distance.
This uncertainty regarding which metric better in a specific scenario further  illustrates the benefits of  predicting the right metric to use in each specific situation.

\begin{figure}
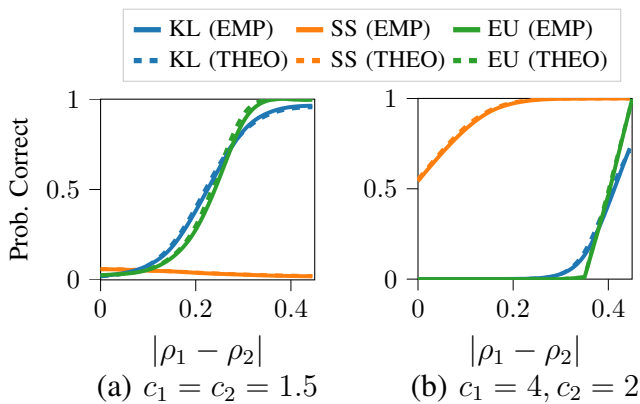

\centering
\begin{subfigure}
    {
\begin{tikzpicture}

\definecolor{darkgray176}{RGB}{176,176,176}
\definecolor{green}{RGB}{0,128,0}
\definecolor{magenta}{RGB}{255,0,255}

\definecolor{subspaceorange}{RGB}{255,127,14}
\definecolor{klblue}{RGB}{31,119,180}
\definecolor{euclideangreen}{RGB}{44,160,44}
\definecolor{lightgray204}{RGB}{204,204,204}

\begin{axis}[
width=0,
height=0,
at={(0,0)},
scale only axis,
scale only axis,
legend style={legend cell align=left, align=left, draw=lightgray204, at={(1.85in,0)},
anchor=south,legend columns=3,
font=\small,
legend image post style={scale=0.6}
}
]
\addplot [semithick, klblue,line width=2pt]
table {%
0 0
};
\addlegendentry{KL (EMP)}
\addplot [semithick, subspaceorange,line width=2pt]
table {%
0 0
};
\addlegendentry{SS (EMP)}
\addplot [semithick, euclideangreen, line width=2pt]
table {%
0 0
};
\addlegendentry{EU (EMP)}

\addplot [semithick, klblue, dashed,line width=2pt]
table {%
0 0
};
\addlegendentry{KL (THEO)}
\addplot [semithick, subspaceorange, dashed, line width=2pt]
table {%
0 0
};
\addlegendentry{SS (THEO)}

\addplot [semithick, euclideangreen, dashed,line width=2pt]
table {%
0 0
};
\addlegendentry{EU (THEO)}
\end{axis}
\end{tikzpicture}}
\vspace{-1em}
\end{subfigure}
\\
\begin{subfigure}
    {\input{figs/final/2final_a}}
\end{subfigure}
\begin{subfigure}
    {\input{figs/final/2final_b}}
\end{subfigure}
\vspace{-1.5em}
\caption{Empirical (EMP) and theoretical (THEO) probability of correct clustering four SCMs into two groups for $M=50$ and $\rho_1 = 0.4$.}
\label{fig:results:clustering_SCMs_xaxis_delta}
\vspace{-1em}
\end{figure}

Our results also allow us to study the separability of the different clusters with the considered functionals, depending on how different the true covariances are. 
Increasing the distance between these parameters can be understood as further increasing the difference between the covariance matrices generating data. To ease interpretation, we consider here a simplified scenario with only four sample covariance matrices, so that $\mathbf{R}_1=\mathbf{R}_2$ and $\mathbf{R}_4=\mathbf{R}_4$. Specifically, Fig.~\ref{fig:results:clustering_SCMs_xaxis_delta} portrays the empirical probability of correct clustering (solid lines) against the theoretical probability (dashed lines) for increasing $\Delta\rho = |\rho_1 - \rho_3|$ and fixed $M = 50, \rho_1 = 0.4$ in two undersampled regimes. 
Even in this simple scenario with four SCMs and two groups, there is no single predominant metric that consistently outperforms the others. Depending on the specific scenario ($\rho_j, c_j, j=1,\ldots, 4$),  different metrics might perform better or worse than another.  In such scenarios, where different metrics will lead to different results, it becomes particularly useful to (theoretically) assess different  algorithms.

\vspace{-0.7em}
\section{Conclusion}

The asymptotic characterization of a general class of functionals involving two SCMs that can be expressed as the sum of traces of analytic functions applied to each matrix separately has been studied. Results are established in the asymptotic regime where both the sample size and the observation dimension tend to infinity at the same rate. These generic results are valid both in the undersampled and oversampled regimes as well as for complex- and real-valued observations. Moreover, we have particularized these results to three commonly used distances between covariance matrices, namely the Euclidean distance, a symmetrized version of the Kullback–Leibler divergence and Subspace Similarity based on the principal angles of the compared subspaces. Numerical simulations confirm the validity of the presented results by illustrating the accuracy of the asymptotic behavior of these metrics when compared to their empirical distributions. Furthermore, we have also illustrated how these results can assist on properly designing clustering solutions. Specifically, we show that CLTs derived in this work become useful when assessing the quality of a clustering for different distance metrics. Finally, we stress that the results presented here are also applicable to describe the behavior of other algorithms that are based on pairwise similarity comparisons and to other distances that follow the same structure as the class of distances discussed throughout this work.


\appendices

\vspace{-1em}
\section{Some useful lemmas} \label{sec:lemmasbounds}
The following lemmas can be proven using conventional tools in random matrix theory. Their proofs are therefore omitted due to space constraints

\begin{lemma}
\label{lemma:UnifBounds}Let $j\in\{0,1\}$. Under \textbf{(As1)}-\textbf{(As3)} we have
\begin{gather}
\sup_{M}\sup_{z\in C_j}\left\Vert \mathbf{Q}_{j}(\omega_{j}\left(  z\right)
)\right\Vert <+\infty\nonumber\\
0<\inf_{M}\inf_{z\in C_j}\left\vert \omega_{j}\left(  z\right)  \right\vert
\leq\sup_{M}\sup_{z\in C_j}\left\vert \omega_{j}\left(  z\right)  \right\vert
<+\infty\label{eq:boundednessw(z)}\\
\inf_{M}\mathrm{dist}\left\{  \omega_{j}\left(  z\right)  ,(-\infty,\mu_{\inf
}^{(j)}]\right\}  >0\nonumber
\end{gather}
where $\mu_{\inf}^{(j)} = \inf_{M}\mu_0^{(j)}$. Furthermore,
\begin{gather}
\sup_{M\geq M_{0}}\sup_{z\in C_j}\left\Vert \mathbf{\hat{Q}}_{j}(z)\right\Vert
<+\infty\nonumber\\
0<\inf_{M\geq M_{0}}\inf_{z\in C_j}\left\vert \hat{\omega}_{j}\left(  z\right)
\right\vert \leq\sup_{M\geq M_{0}}\sup_{z\in C_j}\left\vert \hat{\omega}
_{j}\left(  z\right)  \right\vert <+\infty\label{eq:boundednessw(z)prime}\\
\inf_{M\geq M_{0}}\mathrm{dist}\left\{  \hat{\omega}_{j}\left(  z\right)
,(-\infty,\mu_{\inf}^{(j)}]\right\}  >0\nonumber
\end{gather}
with probability one for some $M_{0}$ sufficiently high.
\end{lemma}

\begin{lemma}
\label{lemma:UnifBounds2}Let $j\in\{0,1\}$. Under \textbf{(As1)}-\textbf{(As3)} we have
\[
\sup_{M}\sup_{z\in C_j} \overline{\Gamma}_j(\omega_j(z)) <1
\]
\end{lemma}

\section{Derivation of the asymptotic second-order mean and variances of euclidean distance}
\label{sec:derivation2ndOrderMeanVar}


Let us evaluate the asymptotic mean and variance of the Euclidean distance between SCMs. To that effect, one must carry out the
integrals in (\ref{eq:asymptMean})-(\ref{eq:asymptVariance}) for
$g(z_{1},z_{2})=\left(  z_{1}-z_{2}\right)  ^{2}$ under the assumption that
all the contours enclose $\{0\}$. Hence, one can follow exactly the
same approach that was used in Section \ref{sec:Euclideandist1stOrder}. The
main idea is to first use the change of variable $z\mapsto\omega=\omega
_{j}(z)$. The resulting contour $\cz{\omega_j}=\omega
_{j}(\cz{j})$ encloses all the singularities of the integrand (see further Figure \ref{fig:contours}, so
that one can apply a second change of variables $\omega\mapsto\zeta
(\omega)=\omega^{-1}\,$\ in a way that $\zeta(\cz{\omega_{j}})$
encloses zero and no other singularity. By direct application of this
technique we may find
\begin{equation*}
\frac{1}{2\pi\mathrm{j}}\oint\nolimits_{\cz{1}}\frac{\omega_{1}
}{z_{1}}\mathfrak{m}\left(  \omega_{1},\omega_{2}\right)  dz_{1}
=
\frac{\mathrm{tr}\left[  \mathbf{\mathbf{R}}_{2}^{2}\mathbf{Q}_{2}^{3}\left(  \omega_{2}\right)  \Omega_{2}\left(  \omega_{2};\mathbf{I}_{M}\right)  \right]}{N_{2}\left(1-\Gamma_{2}(\omega_{2})\right)}
\end{equation*}
together with
\[
\frac{1}{2\pi\mathrm{j}}\oint\nolimits_{\cz{1}}\omega_{1} 
\mathfrak{m}\left(  \omega_{1},\omega_{2}\right)
dz_{1}=\frac{1}{N_{2}}\frac{\mathrm{tr}\left[  \mathbf{\mathbf{R}}_{2}
^{2}\mathbf{Q}_{2}^{3}\left(  \omega_{2}\right)  \Omega_{2}\left(  \omega
_{2};\mathbf{R}_{1}\right)  \right]  }{1-\Gamma_{2}(\omega_{2})}.
\]
Repeating the approach with respect to the second variable $z_{2}$ we find
\[
\frac{1}{\left(  2\pi\mathrm{j}\right)  ^{2}}\oint\nolimits_{\cz{1}}\oint\nolimits_{\cz{2}}z_{1}z_{2}\frac{\omega_{1}
\omega_{2}}{z_{1}z_{2}}\mathfrak{m}\left(  \omega_{1},\omega_{2}\right)
dz_{1}dz_{2}=0
\]
and
\[
\frac{1}{\left(  2\pi\mathrm{j}\right)  ^{2}}\oint\nolimits_{\cz{1}} \oint\nolimits_{\cz{2}}z_{2}^{2}\frac{\omega_{1}
\omega_{2}}{z_{1}z_{2}}\mathfrak{m}\left(  \omega_{1},\omega_{2}\right)
dz_{1}dz_{2}=\frac{1}{N_{2}}\mathrm{tr}\left[  \mathbf{\mathbf{R}}_{2}
^{2}\right].
\]
We can therefore conclude that $\mathfrak{m}_{M}$ takes the form in (\ref{eq:secondOrderMeanEuclidean}).

The same approach can be followed to evaluate the asymptotic variance. In this
case, we will use the fact that for any $\omega^{\prime}\in \cz{\omega_j}=\omega_{j}(\cz{j})$, the function $\Gamma
_{j}(\omega,\omega^{\prime})=1$ has all its solutions inside $\tilde{\cz{j}}$, see~\cite{Schenck22}. Therefore, all the singularities of
$\sigma_{j}^{2}\left(  \omega,\omega^{\prime};\mathbf{A},\mathbf{B}\right)  $
are located inside $\mathcal{C}_{\omega_{j}}^{(l)}$. We can therefore use the
same integration technique as before by first applying the change of variables
$z_{j}\mapsto\omega_{j}=\omega_{j}(z_{j})$ and noting that the resulting
contour $\cz{\omega_j}=\omega_{j}(\cz{j})$
encloses all the singularities of the integrand, because $\cz{j}$
is built to enclose zero. Applying then a second change of variables
$\omega_{j}\mapsto\zeta(\omega_{j})=\omega_{j}^{-1}$ one can see that
$\zeta(\cz{\omega_{j}})$ encloses zero and no other singularity.
Computing the corresponding residue, one easily finds that
\begin{multline*}
\frac{1}{2\pi\mathrm{j}}\oint\nolimits_{\cz{1}} \left(
z_{1}-z_{2}\right)  ^{2}\frac{\omega_{1}}{z_{1}}\sigma_{1}^{2}\left(
\omega_{1},\omega_{1}^{\prime};\mathbf{A},\mathbf{B}\right)  dz_{1}=\\
=\left(  2z_{2}+\omega_{1}^{\prime}\frac{1}{N_{1}}\mathrm{tr}\left[
\mathbf{R}_{1}\mathbf{Q}_{1}\left(  \omega_{1}^{\prime}\right)  \right]
\right) \times \\ \times \frac{1}{N_{1}}\mathrm{tr}\left[  \mathbf{R}_{1}\mathbf{Q}_{1}\left(
\omega_{1}^{\prime}\right)  \mathbf{AR}_{1}\mathbf{Q}_{1}\left(  \omega
_{1}^{\prime}\right)  \Omega_{1}\left(  \omega_{1}^{\prime};\mathbf{B}\right)
\right]  \\
-\frac{1}{N_{1}}\mathrm{tr}\left[  \mathbf{R}_{1}^{2}\mathbf{Q}_{1}\left(
\omega_{1}^{\prime}\right)  \mathbf{AR}_{1}\mathbf{Q}_{1}\left(  \omega
_{1}^{\prime}\right)  \Omega_{1}\left(  \omega_{1}^{\prime};\mathbf{B}\right)
\right]  \\
-\frac{1}{N_{1}}\mathrm{tr}\left[  \mathbf{R}_{1}\mathbf{Q}_{1}\left(
\omega_{1}^{\prime}\right)  \mathbf{AR}_{1}^{2}\mathbf{Q}_{1}\left(
\omega_{1}^{\prime}\right)  \Omega_{1}\left(  \omega_{1}^{\prime}%
;\mathbf{B}\right)  \right]  \\
+\omega_{1}^{\prime}\frac{1}{N_{1}}\mathrm{tr}\left[  \mathbf{R}_{1}%
\mathbf{Q}_{1}\left(  \omega_{1}^{\prime}\right)  \mathbf{A}\right]  \frac
{1}{N_{1}}\mathrm{tr}\left[  \mathbf{R}_{1}^{2}\mathbf{Q}_{1}^{2}\left(
\omega_{1}^{\prime}\right)  \Omega_{1}\left(  \omega_{1}^{\prime}%
;\mathbf{B}\right)  \right]  .
\end{multline*}
We can now multiply by $\left(  z_{1}^{\prime}-z_{2}^{\prime}\right)
^{2}\omega_{1}^{\prime}/z_{1}^{\prime}$ and integrate with respect to
$z_{1}^{\prime}$. Following exactly the same approach as above we obtain,
after some algebra, the integral in (\ref{eq:integralABsquares2}) at the top of the next page. 

Now, to solve the integrals with respect to $z_2$, we apply the following result
\begin{align*}
\frac{1}{2\pi\mathrm{j}}\oint\nolimits_{\cz{2}}\mathbf{Q}
_{2}(\omega_{2})\frac{\omega_{2}}{z_{2}}dz_{2} &  =\mathbf{I}_{M}
\\
\frac{1}{2\pi\mathrm{j}}\oint\nolimits_{\cz{2}} z_{2}
\mathbf{Q}_{2}(\omega_{2})\frac{\omega_{2}}{z_{2}}dz_{2} &  =\mathbf{R}_{2}.
\end{align*}
A direct application of the above identities allows us to obtain the integral
in (\ref{eq:finaldoublepowerint}) at the top of the next page. 
Obviously, the integral of the term $\sigma_{2}^{2}\left(  \omega_{2}%
,\omega_{2}^{\prime};\mathbf{Q}_{1}(\omega_{1}),\mathbf{Q}_{1}(\omega
_{1}^{\prime})\right)  $ can be obtained from the above by simply swapping the
two indices.

We finally evaluate the last integral, first by noting that, using the same
integration technique, we obtain%
\begin{multline*}
\frac{1}{2\pi\mathrm{j}}\oint\nolimits_{\cz{1}}\left(
z_{1}-z_{2}\right)  ^{2}\frac{\omega_{1}}{z_{1}}\frac{\mathrm{tr}^{2}\left[
\mathbf{AQ}_{1}\left(  \omega_{1}\right)  \right]  }{1-\Gamma_{1}(\omega
_{1},\omega_{1}^{\prime})}dz_{1}=\\
=\left(  2z_{2}-\frac{1}{N_{1}}\mathrm{tr}\left[  \mathbf{R}_{1}\right]
\right)  \mathrm{tr}^{2}\left[  \mathbf{A}\right]  -2\mathrm{tr}\left[
\mathbf{A}\right]  \mathrm{tr}\left[  \mathbf{AR}_{1}\right]  \\
+\mathrm{tr}^{2}\left[  \mathbf{A}\right]  \frac{1}{N_{1}}\mathrm{tr}\left[
\mathbf{R}_{1}^{2}\mathbf{Q}_{1}\left(  \omega_{1}^{\prime}\right)  \right]
\end{multline*}
and therefore replacing $\mathbf{A}$\ with $\mathbf{BQ}_{2}\left(  \omega
_{2}\right)  \mathbf{A}$ and using the same technique we find%
\begin{multline*}
\frac{1}{\left(  2\pi\mathrm{j}\right)  ^{2}}\oint\nolimits_{\cz{1}}\oint\nolimits_{\cz{2}} \left(  z_{1}-z_{2}\right)
^{2}\frac{\omega_{1}}{z_{1}}\frac{\omega_{2}}{z_{2}}\times\\
\times\frac{\mathrm{tr}^{2}\left[  \mathbf{AQ}_{1}\left(  \omega_{1}\right)
\mathbf{BQ}_{2}\left(  \omega_{2}\right)  \right]  }{\left(  1-\Gamma
_{1}(\omega_{1},\omega_{1}^{\prime})\right)  \left(  1-\Gamma_{2}(\omega
_{2},\omega_{2}^{\prime})\right)  }dz_{1}dz_{2}=-2\mathrm{tr}^{2}\left[
\mathbf{AB}\right]
\end{multline*}
and consequently, after some manipulation,
\begin{multline*}
  \frac{1}{\left(  2\pi\mathrm{j}\right)  ^{4}}\oint\nolimits_{\cz{1}}\oint\nolimits_{\cz{2}}\oint\nolimits_{\cz{1}}\oint\nolimits_{\cz{2}}\left(  z_{1}-z_{2}\right)
^{2}\left(  z_{1}^{\prime}-z_{2}^{\prime}\right)  ^{2}\frac{\omega_{1}\omega_{1}^{\prime} \omega_{2}\omega_{2}^{\prime} }{z_{1} z_{1}^{\prime} z_2 z_{2}^{\prime}} \times\\
 \times
\varrho (\omega_1,\omega_1^{\prime},\omega_2,\omega_2^{\prime}) 
dz_{1}dz_{2}dz_{1}^{\prime}dz_{2}^{\prime} =4\frac{\mathrm{tr}^{2}[\mathbf{R}_{1}\mathbf{R}_{2}]}{N_{1}N_{2}}.
\end{multline*}
Adding the three integrals, we can conclude that the asymptotic variance
takes the expression in (\ref{eq:varianceEuclidean}).

 \begin{figure*}[t]
     \centering
   \begin{multline}
  \frac{1}{\left(  2\pi\mathrm{j}\right)  ^{2}}\oint\nolimits_{\cz{1}}\oint\nolimits_{\cz{1}}\left(  z_{1}-z_{2}\right)
^{2}\left(  z_{1}^{\prime}-z_{2}^{\prime}\right)  ^{2}\frac{\omega_{1}}{z_{1}%
}\frac{\omega_{1}^{\prime}}{z_{1}^{\prime}}\sigma_{1}^{2}\left(  \omega
_{1},\omega_{1}^{\prime};\mathbf{A},\mathbf{B}\right)  dz_{1}dz_{1}^{\prime}  =\\
  =\frac{1}{N_{1}}\mathrm{tr} \Big[\mathbf{R}_{1}\left(  \left(  \frac{1}{N_{1}%
}\mathrm{tr}\left[  \mathbf{R}_{1}\right]  -2z_{2}\right)  \mathbf{A}%
+\mathbf{R}_{1}\mathbf{A+AR}_{1}+\frac{1}{N_{1}}\mathrm{tr}\left[
\mathbf{R}_{1}\mathbf{A}\right]  \mathbf{I}_{M}\right)  \times\\
  \times\mathbf{R}_{1}\left(  \left(  \frac{1}{N_{1}}\mathrm{tr}\left[
\mathbf{R}_{1}\right]  -2z_{2}^{\prime}\right)  \mathbf{B}+\mathbf{BR}%
_{1}+\mathbf{R}_{1}\mathbf{B+}\frac{1}{N_{1}}\mathrm{tr}\left[  \mathbf{R}%
_{1}\mathbf{B}\right]  \mathbf{I}_{M}\right)  \Big] \\
  +\frac{1}{N_{1}}\mathrm{tr}\left[  \mathbf{R}_{1}^{2}\mathbf{A}\right]
\frac{1}{N_{1}}\mathrm{tr}\left[  \mathbf{R}_{1}^{2}\mathbf{B}\right]
+\frac{1}{N_{1}}\mathrm{tr}\left[  \mathbf{R}_{1}^{2}\right]  \frac{1}{N_{1}%
}\mathrm{tr}\left[  \mathbf{R}_{1}\mathbf{AR}_{1}\mathbf{B}\right]   \label{eq:integralABsquares2}
\end{multline}

\begin{multline}
 \frac{1}{\left(  2\pi\mathrm{j}\right)  ^{4}}\oint\nolimits_{\cz{1}}\oint\nolimits_{\cz{1}}\oint\nolimits_{\cz{2}}\oint\nolimits_{\cz{2}}
 \left(  z_{1}-z_{2}\right)
^{2}\left(  z_{1}^{\prime}-z_{2}^{\prime}\right)  ^{2}\left(  \frac{\omega
_{1}}{z_{1}}\frac{\omega_{1}^{\prime}}{z_{1}^{\prime}}\frac{\omega_{2}}{z_{2}
}\frac{\omega_{2}^{\prime}}{z_{2}^{\prime}}\right) \sigma_{1}^{2}\left(
\omega_{1},\omega_{1}^{\prime};\mathbf{A},\mathbf{B}\right)  dz_{1}
dz_{1}^{\prime}dz_{2}dz_{2}^{\prime}\\
  =2\left(  \frac{1}{N_{1}}\mathrm{tr}\left[  \mathbf{R}_{1}^{2}\right]
\right)  ^{2} +4\frac{1}{N_{1}}\mathrm{tr}\left[  \mathbf{R}_{1}\left(
\mathbf{R}_{1}+\frac{1}{N_{1}}\mathrm{tr}\left[  \mathbf{R}_{1}\right]
\mathbf{I}_{M}-\mathbf{R}_{2}\right)  \mathbf{R}_{1}\left(  \mathbf{R}
_{1}+\frac{1}{N_{1}}\mathrm{tr}\left[  \mathbf{R}_{1}\right]  \mathbf{I}
_{M}-\mathbf{R}_{2}\right)  \right] \label{eq:finaldoublepowerint}
\end{multline}
     \rule[2ex]{\linewidth}{1pt} 
 \end{figure*}

\bibliographystyle{IEEEtran}
\bibliography{./bib/IEEEbib}

\end{document}